\documentclass[10pt,a4paper]{amsart}
\usepackage{amssymb,amsmath}
\usepackage[american,french]{babel}
\usepackage{amsopn}
\usepackage{graphicx}
\usepackage[ec]{aeguill}
\usepackage{fancyhdr}
\usepackage[T1]{fontenc}
\title{On the pseudospectrum of elliptic quadratic differential operators}
\newcommand{\rr}{\mathbb{R}}
\newcommand{\eps}{\varepsilon}
\newcommand{\nn}{\mathbb{N}}
\newcommand{\cc}{\mathbb{C}}

\newcommand{\lde}{L^2(\rr^n)}

\def\init{\setcounter{equa}{0}}
\def\inc{\stepcounter{equa}}
\def\num{\tag{\thesubsection.\theequa}}
\begin{document}
\newcounter{equa}
\selectlanguage{american}
\begin{center}
{\large\textbf{ON THE PSEUDOSPECTRUM OF ELLIPTIC QUADRATIC DIFFERENTIAL OPERATORS}\\
\bigskip
\medskip
Karel \textsc{Pravda-Starov}\\
\bigskip
University of California, Berkeley}
\end{center}
\bigskip
\bigskip

\newtheorem{lemma}{Lemma}[subsection]
\newtheorem{definition}{Definition}[subsection]
\newtheorem{proposition}{Proposition}[subsection]
\newtheorem{theorem}{Theorem}[subsection]

\textbf{Abstract}. We study the pseudospectrum of a class of non-selfadjoint differential operators. Our work consists in a detailed study of the microlocal properties, which rule the spectral 
stability or instability phenomena appearing under small perturbations for elliptic quadratic differential operators. The class of elliptic quadratic differential operators stands for the class of 
operators defined in the Weyl quantization by
complex-valued elliptic quadratic symbols. We establish in this paper a simple necessary and sufficient condition on the Weyl symbol of these operators, which ensures the stability of their spectra. 
When this condition is violated, we prove that it occurs some strong spectral instabilities for the high energies of these operators, in some regions which can be far away from their spectra. 
We give a precise geometrical description of them, which explains the results obtained for these operators in some numerical simulations giving the computation of \og false eigenvalues \fg \ 
far from their spectra by algorithms for eigenvalues computing.

\medskip


\medskip


\section{Introduction}
 
\subsection{Miscellaneous facts about pseudospectrum}
\init

In recent years, there has been a lot of interest in studying the pseudospectrum of non-selfadjoint operators. The study of this notion has
been initiated by noticing that for certain problems of science and engineering involving non-selfadjoint operators, the predictions suggested
by spectral analysis do not match with the numerical simulations. This fact lets thinking that in some cases the only knowledge of the spectrum of an operator is not enough to understand 
sufficiently its action. To supplement this lack of information contained in the spectrum, some new subsets of the complex plane called pseudospectra have been defined. The main idea 
about the definition of these new subsets is that it is interesting to study not only the points where the resolvent of an operator is not defined, i.e. its spectrum, but also where this resolvent 
is large in norm. This explains the following definition of the $\varepsilon$-pseudospectrum $\sigma_{\varepsilon}(A)$ of a matrix or an operator $A$, 
$$\sigma_{\varepsilon}(A)=\Big\{z \in \cc,  \ \|(zI-A)^{-1}\| \geq \frac{1}{\varepsilon} \Big\},$$
for any $\eps>0$, if we write by convention that $\|(zI-A)^{-1}\|=+\infty$ for every point $z$ belonging to the spectrum $\sigma(A)$ of the operator. 

Let us mention that there exists an abundant literature about this notion of pseudospectrum. We refer here for the definition and some general properties of pseudospectra to the paper
\cite{trefethen} of L.N. Trefethen. Let us also point out the more recently published book \cite{trefethen2}, which draws up a wide all-round view of this topic and gives a lot of illustrations.

According to the previous definition, studying the pseudospectra of an operator is exactly studying the level lines of the norm of its resolvent. What is interesting in studying such level lines is that it 
gives some information about the spectral stability of the operator. Indeed, pseudospectra can be defined in an equivalent way in term of spectra of perturbations of the operator. For instance, we have for 
any $A \in M_n(\cc)$,
$$\sigma_{\varepsilon}(A)=\{z \in \cc, \ z \in \sigma(A+B) \ \textrm{for some} \ B \in M_n(\cc) \ \textrm{with} \ \|B\| \leq \varepsilon \}.$$
It follows that a complex number $z$ belongs to the $\varepsilon$-pseudospectrum of a matrix $A$ if and only if it belongs to the spectrum of one of its perturbations $A+B$ with $\|B\| \leq \varepsilon$. More generally, if $A$ is a closed unbounded linear operator with a dense domain on a complex Hilbert space $H$, the result of Roch and Silbermann in 
\cite{roch}  gives that 
$$\sigma_{\eps}(A)=\bigcup_{B \in \mathcal{L}(H), \ \|B\|_{\mathcal{L}(H)} \leq \eps}{\sigma(A+B)},$$
where $\mathcal{L}(H)$ stands for the set of bounded linear operators on $H$. From this second description, we understand the interest in studying such subsets if we
want for example to compute numerically some eigenvalues of an operator. Indeed, we start to do it by discretizing this operator. This discretization and inevitable round-off errors will generate 
some perturbations of the initial operator. Eventually, algorithms for eigenvalues computing will determine the eigenvalues of a perturbation of the initial operator, i.e. a value in 
a $\varepsilon$-pseudospectrum of the initial operator but not necessarily a spectral one. This explains why it is important in such numerical computations to understand 
if the $\eps$-pseudospectra of studied operators contain more or less deeply their spectra.

Let us first notice that this study is a priori non-trivial only for non-selfadjoint operators, or more precisely for \textit{non-normal} operators. Indeed, we have for a \textit{normal} operator $A$ 
an exact expression of the norm of its resolvent given by the following classical formula 
(see for example (V.3.31) in \cite{kato}),
\begin{equation}\label{1}\inc
\forall z \not\in \sigma(A), \ \|(zI-A)^{-1}\| = \frac{1}{d\big(z,\sigma(A)\big)}, \num
\end{equation}
where $d\big(z,\sigma(A)\big)$ stands for the distance between $z$  and the spectrum of the operator, when $A$ is a closed unbounded linear 
operator with a dense domain on a complex Hilbert space. This formula proves that the resolvent of a \textit{normal} operator cannot 
blow up far from its spectrum. It ensures the stability of its spectrum under small perturbations because the $\eps$-pseudospectrum 
is exactly equal in this case to the $\eps$-neighbourhood of the spectrum
\begin{equation}\label{2}\inc
\sigma_{\eps}(A) = \big\{z \in \cc :  d\big(z,\sigma(A) \big) \leq \eps \big\}. \num
\end{equation}
Nevertheless it is well-known that this formula (\ref{1}) is no more true for non-normal operators. For such operators, it can occur that their resolvents are very large in norm 
far from their spectra. This induces that the spectra of these operators can be very unstable under small perturbations. To illustrate this fact, let us consider the 
case of the rotated harmonic oscillator and the following 
numerical computation of its spectrum. The rotated harmonic oscillator is a simple example of elliptic quadratic differential operator
$$H_{c}=D_x^2+c x^2, \ D_x=i^{-1}\partial_x,$$
with $c=e^{i \pi/4}$. The numerical computation is performed on the matrix discretization
$$\big{(}(H_{c} \Psi_i,\Psi_j)_{L^2(\rr)}\big{)}_{1 \leq i,j \leq N},$$ 
where $N$ is an integer taken equal to $100$ and $(\Psi_j)_{j \in \mathbb{N}^*}$ stands for the basis of $L^{2}(\rr)$ composed by Hermite functions. 
\begin{figure}[ht]
\caption{Computation of some level lines of the norm of the resolvent $\|(H_c-z)^{-1}\|=\eps^{-1}$ for the rotated harmonic oscillator $H_c$ with $c=e^{i\pi/4}$. The right column gives the corresponding values of 
$\log_{10} \eps$.}
\centerline{\includegraphics[scale=0.6]{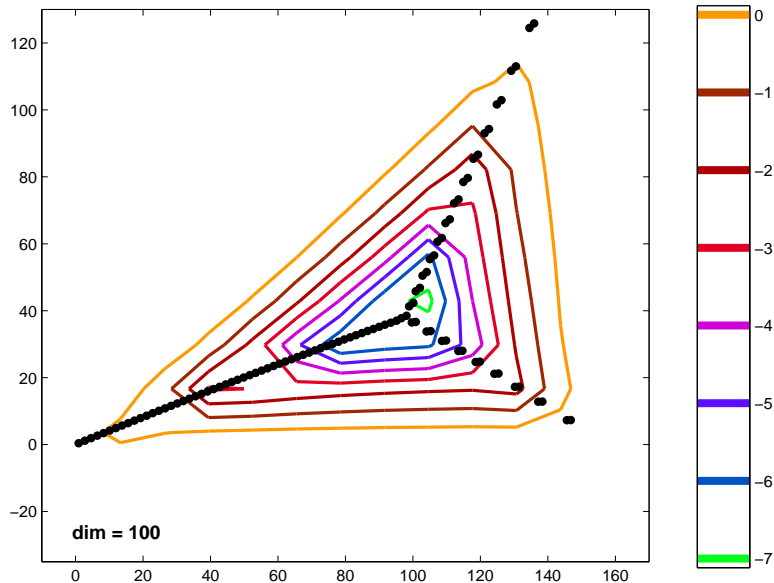}}
\end{figure}  
The black dots appearing on this computation stand for the numerically computed eigenvalues. We can notice on this numerical simulation that the computed low energies are very close to theoretical ones
since the spectrum of the rotated harmonic oscillator is only composed of eigenvalues regularly spaced out on the half-line $e^{i \pi/8} \rr_+^*$,
$$\sigma(H_c)=\{e^{i \pi/8}(2n+1) :  n \in \nn\}.$$  
However we notice that it is no more true for the high energies. It occurs for them some strong spectral instabilities, which lead to the computation of \og false eigenvalues \fg \ 
far from the half-line $e^{i \pi/8}\rr_+^*$. Let us mention that some comparable computations can be found in \cite{daviesosc}. In this paper, we are interested in studying when and how this kind of phenomena
occurs in the class of elliptic quadratic differential operators.

\subsection{Elliptic quadratic differential operators}
\init

We study here the class of elliptic quadratic differential operators. It is the class of pseudodifferential operators defined in the Weyl quantization
\begin{equation}\label{3}\inc
q(x,\xi)^w u(x) =\frac{1}{(2\pi)^n}\int_{\rr^{2n}}{e^{i(x-y).\xi}q\Big(\frac{x+y}{2},\xi\Big)u(y)dyd\xi}, \num
\end{equation}
by some symbols $q(x,\xi)$, where $(x,\xi) \in \rr^{n} \times \rr^n$ and $n \in \nn^*$, which are some \textit{complex-valued elliptic quadratic forms} i.e. complex-valued quadratic forms verifying
\begin{equation}\label{3.5}\inc
(x,\xi) \in \rr^{n} \times \rr^n, \ q(x,\xi)=0 \Rightarrow (x,\xi)=(0,0). \num
\end{equation}
Let us first notice that since the symbols of these operators are some quadratic forms, these are only some differential operators, which are a priori non-selfadjoint because their Weyl
symbols are complex-valued. As mentioned before, the rotated harmonic oscillator is an example of such an operator since we have 
$$D_x^2+e^{i\theta} x^2=(\xi^2+e^{i\theta}x^2)^w, \ 0<\theta<\pi,$$
if $D_x=i^{-1}\partial_x$. This operator is a very simple example of non-selfadjoint operator for which we have noticed on the previous numerical simulation that 
it occurs some strong spectral instabilities under small perturbations for its high energies. These phenomena have been studied in several recent works. We can mention in particular the 
works of L.S. Boulton \cite{boulton}, E.B. Davies \cite{daviesosc}, K. Pravda-Starov \cite{karel3} and M. Zworski~\cite{zworski2}, which have given a good understanding of these phenomena.

A question, which has been at the origin of this work, has been to study if these phenomena peculiar to the rotated harmonic oscillator are representative, or not, of what occurs more
generally in the class of elliptic quadratic differential operators in every dimension. We have tried to answer to the following questions:
\begin{itemize}
\item[-] Does it always occur some strong spectral instabilities under small perturbations for the high energies of these operators ?
\item[-] If it is not the case, is it possible to give a necessary and sufficient condition on the Weyl symbols of these operators, which ensures their spectral stability~?
\item[-] Can we precisely describe the geometry, which separates the regions of the resolvent sets where the resolvents of these operators blow up in norm from the ones where one keeps
a control on their sizes ? 
\end{itemize}

To understand these spectral stability or instability phenomena, we need to study the microlocal properties, which rule these phenomena in the class of elliptic quadratic differential operators.
Let us mention that it is M. Zworski who first underlined in~\cite{zworski2} the close link between these questions of spectral instabilities and some results of 
microlocal analysis about the solvability of pseudodifferential operators.

\subsection{Semiclassical pseudospectrum}
\init

To answer to these previous questions, it is interesting to use a semiclassical setting and to study a notion of pseudospectrum in this new setting. We define for 
a semiclassical family $(P_h)_{0<h \leq 1}$ of operators on $L^2(\rr^n)$, with a domain~$D$, the following notions of semiclassical pseudospectra. 

\begin{definition}\label{4}
For all $\mu \geq 0$, the set  
$$
\Lambda_{\mu}^{\emph{\textrm{sc}}}
(P_h)=\big\{z \in \mathbb{C} : \forall C >0, \forall h_0>0, \exists \ 0<h<h_0,
\ \|(P_h-z)^{-1}\| \geq C h^{-\mu}  \big\},
$$
is called \emph{semiclassical pseudospectrum of index $\mu$} of the semiclassical family $(P_h)_{0<h \leq 1}$.  
The \emph{semiclassical pseudospectrum of infinite index} is defined by 
$$\Lambda_{\infty}^{\emph{\textrm{sc}}}(P_h)=\bigcap_{\mu \geq 0}{\Lambda_{\mu}^{\emph{\textrm{sc}}}(P_h)}.$$
\end{definition}

\bigskip

\noindent
With this definition, the points in the complement of the semiclassical pseudospectrum of index $\mu$ are the points of the complex plane
where we have the following control of the resolvent's norm for sufficiently small values of the semiclassical parameter $h$,
\begin{equation}\label{4.5}\inc
\exists C>0, \exists h_0>0, \forall \ 0<h<h_0, \ \|(P_h-z)^{-1}\|<C h^{-\mu}. \num
\end{equation}
To prove the existence of semiclassical pseudospectrum of index $\mu$, we will study the question of existence of semiclassical quasimodes 
\inc
\begin{multline}\label{5}
\forall C>0, \forall h_0>0, \exists \ 0<h<h_0, \exists u_h \in D, \\
 \|u_h\|_{L^2(\rr^n)}=1 \ \textrm{and }  \|P_h u_h-zu_h\|_{L^2(\rr^n)} \leq C h^{\mu}, \num
\end{multline}
in some points $z$ of the resolvent set, which can be considered as some \og almost eigenvalues \fg \ in $O(h^{\mu})$ in the semiclassical limit.
Let us notice that the definition chosen here for the notions of semiclassical pseudospectra differ from the one given in \cite{dencker} for a semiclassical pseudodifferential
operator. In fact, we have chosen a definition for semiclassical pseudospectra inspired by the remark made p.388 in \cite{dencker}, because this definition only depends 
on the properties of the semiclassical operator rather than on its symbol.

The interest of working in a semiclassical setting is a matter of geometry. We can explain this choice by the fact that it is easier for an elliptic quadratic differential 
operator $q(x,\xi)^w$ to describe the geometry of semiclassical pseudospectra of its associated semiclassical operator 
$(q(x,h\xi)^w)_{0<h \leq 1},$ 
than to describe directly the geometry of its $\eps$-pseudospectra. The semiclassical setting is particularly well-adapted for the study of elliptic 
quadratic differential operators because there exists a simple link between this semiclassical setting and the quantum one. Indeed, using that the symbols of these
operators are some quadratic forms $q$, we obtain from the change of variables, $y=h^{1/2}x$ with $h>0$, the following identity between the quantum operator 
$q(x,\xi)^w$ and its associated semiclassical operator $(q(x,h\xi)^w)_{0<h \leq 1}$,
\begin{equation}\label{6}\inc 
q(x,\xi)^w-\frac{z}{h}=\frac{1}{h}\big{(}q(y,h\eta)^w-z \big{)}, \num
\end{equation}
if $z \in \cc$. This identity allows to get some information about the resolvent's norm behaviour of the quantum operator  
$$\big(q(x,\xi)^w-z\big)^{-1},$$ 
if we have some information about semiclassical pseudospectra for its associated semiclassical operator. Let us mention for example that if a non-zero complex number $z$
belongs to the semiclassical pseudospectrum of infinite index of the operator 
$$(q(x,h\xi)^w)_{0<h \leq 1},$$ 
the identity (\ref{6}) induces that the resolvent's norm of the quantum
operator blows up along the half-line $z \rr_+$ with a rate faster than any polynomials
\begin{equation}\label{7}\inc 
\forall N \in \nn, \forall C>0, \forall \eta_0 \geq 1, \exists \eta \geq \eta_0, \ \|\big(q(x,\xi)^w-z \eta \big)^{-1}\| \geq C \eta^N, \num
\end{equation}
and this, even if this half-line $z\rr_+$ does not intersect the spectrum of the operator~$q(x,\xi)^w$. Conversely, in the case where 
$z \not\in \Lambda_{\mu}^{\textrm{sc}}\big{(}q(y,h\eta)^w\big{)}$, $z \neq 0$ and $0 \leq \mu \leq 1$, 
the identity (\ref{6}) shows that we can find some positive constants $C_1$ and $C_2$ such that the resolvent of the operator $q(x,\xi)^w$ remains bounded in norm in some 
regions of the resolvent set of the shape 
\begin{equation}\label{8}\inc
\big{\{}u \in \cc :  |u| \geq C_1, \ d(\Delta,u) \leq C_2 |\textrm{proj}_{\Delta} u|^{1-\mu} \big{\}} \cap
\ \cc \setminus  \sigma\big{(}q(x,\xi)^w\big{)}, \num
\end{equation} 
where $\Delta=z\rr_+$ and $\textrm{proj}_{\Delta} u$ stands for the orthogonal projection of $u$ on the closed half-line $\Delta$.
Indeed, we obtain from (\ref{4.5}) and (\ref{6}) that 
$$\exists C>0, \exists \eta_0 \geq 1, \forall \eta \geq \eta_0, \ \big\|\big(q(x,\xi)^w-\eta e^{i \textrm{arg}z} \big)^{-1} \big\| < C
\eta^{\mu-1},$$
which induces that for all $v \in D\big(q(x,\xi)^w\big)$ and $\eta \geq \eta_0$,
$$\big\|\big(q(x,\xi)^w-\eta e^{i \textrm{arg}z} \big)v \big\|_{L^2(\rr^n)} \geq C^{-1}\eta^{1-\mu}\|v\|_{L^2(\rr^n)},$$   
if $D\big(q(x,\xi)^w\big)$ stands for the domain of the operator $q(x,\xi)^w$. Then, we can find a constant $\tilde{\eta}_0 \geq 1$ such that if $\tilde{z}$ belongs to 
\begin{multline*}
\big\{u \in \cc : |u| \geq \tilde{\eta}_0, \ d(e^{i \textrm{arg}z}\rr_+,u) \leq
2^{-1}C^{-1}|\textrm{proj}_{e^{i\textrm{arg} z}\rr_+} u|^{1-\mu} \big\} \cap \ \cc \setminus \sigma\big(q(x,\xi)^w \big),
\end{multline*}
then 
$$|\textrm{proj}_{e^{i\textrm{arg} z}\rr_+} \tilde{z}| \geq \eta_0.$$
This induces using the previous estimates and the triangular inequality that if $\tilde{z}$ belongs to 
\begin{multline*}
\big\{u \in \cc : |u| \geq \tilde{\eta}_0, \ d(e^{i \textrm{arg}z}\rr_+,u) \leq
2^{-1}C^{-1}|\textrm{proj}_{e^{i\textrm{arg} z}\rr_+} u|^{1-\mu} \big\} \cap \ \cc \setminus \sigma\big(q(x,\xi)^w \big),
\end{multline*}
we have for all $v \in D\big(q(x,\xi)^w\big)$,
\begin{eqnarray*}
\big\|\big(q(x,\xi)^w-\tilde{z}\big)v \big\|_{L^2} & \geq & 
\big\|\big(q(x,\xi)^w-\textrm{proj}_{e^{i \textrm{arg}z}\rr_+}\tilde{z} \big)v \big\|_{L^2}-d\big(e^{i
\textrm{arg}z}\rr_+,\tilde{z}\big)\|v\|_{L^2}\\
&\geq & 2^{-1} C^{-1} |\textrm{proj}_{e^{i \textrm{arg}z}\rr_+}\tilde{z}|^{1-\mu} \|v\|_{L^2}\\
& \geq& 2^{-1} C^{-1} \eta_0^{1-\mu} \|v\|_{L^2},
\end{eqnarray*}
because $\mu \leq 1$. This last estimate shows that the resolvent of the operator $q(x,\xi)^w$ is bounded in norm by $2C \eta_0^{\mu-1}$ on the set 
\begin{multline*}
\big\{u \in \cc : |u| \geq \tilde{\eta}_0, \ d(e^{i \textrm{arg}z}\rr_+,u) \leq 2^{-1}C^{-1}|\textrm{proj}_{e^{i\textrm{arg} z}\rr_+} u|^{1-\mu} \big\} \cap \ \cc \setminus \sigma\big(q(x,\xi)^w \big).
\end{multline*}
We notice that depending directly on the value of the index $\mu$, $0 \leq \mu<1$, the previous set contains more or less deeply in its interior the half-line
$$\{u \in \cc : |u| \geq \tilde{\eta}_0, \ u \in z \rr_+\}.$$
This fact explains why in the following we will precise carefully the index of the semiclassical pseudospectrum to which a point does not belong when there is no semiclassical 
pseudospectrum of infinite index in that point.

\section{Statement of the results}

\subsection{Some notations and some preliminary facts about elliptic quadratic differential operators}\label{b1}
\init 

Let us begin by giving some notations and recalling known results about elliptic quadratic differential operators. Let $q$ be a complex-valued elliptic quadratic form
\begin{eqnarray*}
q : \rr_x^n \times \rr_{\xi}^n &\rightarrow& \cc\\
 (x,\xi) & \mapsto & q(x,\xi),
\end{eqnarray*}
with $n \in \nn^*$, i.e. a complex-valued quadratic form verifying (\ref{3.5}). The \textit{numerical range} $\Sigma(q)$ of $q$ is defined by the subset in the 
complex plane of all values taken by this symbol  
\begin{equation}\label{9}\inc
\Sigma(q)=q(\rr_x^n \times \rr_{\xi}^n), \num
\end{equation}
and the \textit{Hamilton map} $F \in M_{2n}(\cc)$ associated to the quadratic form $q$ is uniquely defined by the identity
\begin{equation}\label{10}\inc
q\big{(}(x,\xi);(y,\eta) \big{)}=\sigma \big{(}(x,\xi),F(y,\eta) \big{)}, \ (x,\xi) \in \rr^{2n},  (y,\eta) \in \rr^{2n}, \num
\end{equation}
where $q\big{(}\textrm{\textperiodcentered};\textrm{\textperiodcentered} \big{)}$ stands for the polar form associated to the quadratic form $q$ and $\sigma$ is the symplectic
form on $\rr^{2n}$,
\begin{equation}\label{11}\inc
\sigma \big{(}(x,\xi),(y,\eta) \big{)}=\xi.y-x.\eta, \ (x,\xi) \in \rr^{2n},  (y,\eta) \in \rr^{2n}. \num
\end{equation}
Let us first notice that this Hamilton map $F$ is skew-symmetric with respect to $\sigma$. This is just a consequence of the properties of skew-symmetry of the symplectic form
and symmetry of the polar form
\begin{equation}\label{12}\inc
\forall X,Y \in \rr^{2n}, \ \sigma(X,FY)=q(X;Y)=q(Y;X)=\sigma(Y,FX)=-\sigma(FX,Y).\num
\end{equation}

Under this assumption of ellipticity, the numerical range of a quadratic form can only take some very particular shapes. It is a consequence of the following result proved by 
J. Sjöstrand (Lemma 3.1 in \cite{sjostrand}),

\bigskip

\begin{proposition}\label{13}
Let $q :  \rr_x^n \times \rr_{\xi}^n \rightarrow \cc$ a complex-valued elliptic quadratic form. If $n \geq 2$, then there exists $z \in \cc^*$ such that  
$\emph{\textrm{Re}}(z q)$ is a positive definite quadratic form. If $n=1$, the same result is fulfilled if we assume besides that $\Sigma(q) \neq \cc$. 
\end{proposition}

\bigskip
\noindent
This proposition shows that the numerical range of an elliptic quadratic form can only take two shapes. The first possible shape is when $\Sigma(q)$ is equal to the whole 
complex plane. This case can only occur in dimension $n=1$. The second possible shape is when $\Sigma(q)$ is equal to a closed angular sector with a top in $0$ and an 
opening strictly lower than $\pi$.
\begin{figure}[ht]
\caption{Shape of the numerical range $\Sigma(q)$ when $\Sigma(q) \neq \cc$.}
\centerline{\includegraphics[scale=0.9]{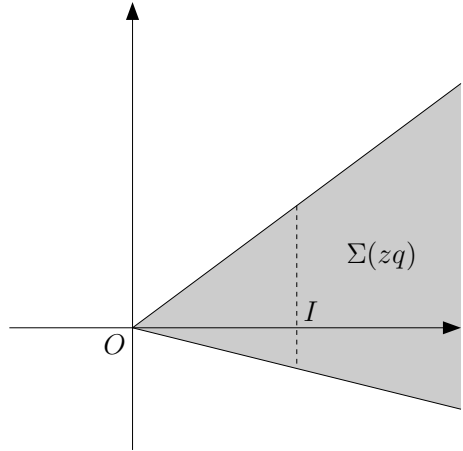}}
\end{figure}

\noindent
Indeed, if $\Sigma(q) \neq \cc$, using that the set $\Sigma(q)$ is a semi-cone  
$$t q(x,\xi)=q(\sqrt{t} x,\sqrt{t}\xi), \ t \in \rr_+, \ (x,\xi) \in \rr^{2n},$$ 
because $q$ is a quadratic form, we have  
$$\Sigma(q)=\rr_+ z^{-1}I,$$
if $z$ is the non-zero complex number given by the proposition \ref{13} and $I$ is the compact interval 
$$I=1+i \ \textrm{Im}(z q)(K),$$
where $K$ is the following compact subset of $\rr^{2n}$, 
$$\big\{(x,\xi) \in \rr^{2n} : \textrm{Re}(z q)(x,\xi)=1\big\}.$$
The compactness of $K$ is a direct consequence of the fact that $\textrm{Re}(zq)$ is a positive definite quadratic form.

Elliptic quadratic differential operators define some Fredholm operators (see Lemma~3.1 in \cite{hypoelliptic} or Theorem 3.5 in \cite{sjostrand}),
\begin{equation}\label{14}\inc
q(x,\xi)^w +z : B \rightarrow \lde, \num
\end{equation}
where $B$ is the Hilbert space 
\begin{equation}\label{14.1}\inc
\big\{ u \in \lde : x^{\alpha} D_x^{\beta} u \in \lde \ \textrm{if} \ |\alpha+\beta| \leq 2\big\}, \num
\end{equation}
with the norm  
$$\|u\|_B^2=\sum_{|\alpha+\beta| \leq 2}{\|x^{\alpha} D_x^{\beta} u\|_{\lde}^2}.$$
The Fredholm index of the operator $q(x,\xi)^w+z$ is independent of $z$ and is equal to~$0$ if $n \geq 2$. In the case where $n=1$, this index can take the values $-2$, $0$ or $2$. More 
precisely, this index is always equal to $0$ if $\Sigma(q) \neq \cc$.

In the following, we will always assume that $\Sigma(q) \neq \cc$. Under this assumption, J.~Sjöstrand has proved in the theorem 3.5 in \cite{sjostrand} (see also Lemma $3.2$ and Theorem~3.3 
in~\cite{hypoelliptic}) that the spectrum of an elliptic quadratic differential operator  
$$q(x,\xi)^w : B \rightarrow \lde,$$ 
is only composed of eigenvalues with finite multiplicity 
\begin{equation}\label{15}\inc
\sigma\big{(}q(x,\xi)^w\big{)}=\Big\{ \sum_{\substack{\lambda \in \sigma(F), \\  -i \lambda \in \Sigma(q) \setminus \{0\}}}
{\big{(}r_{\lambda}+2 k_{\lambda}
\big{)}(-i\lambda) : k_{\lambda} \in \nn} 
\Big\}, \num
\end{equation}
where $F$ is the Hamilton map associated to the quadratic form $q$ and $r_{\lambda}$ is the dimension of the space of generalized eigenvectors of $F$ in $\cc^{2n}$
belonging to the eigenvalue $\lambda \in \cc$. Let us notice that the spectra of these operators is always included in the numerical range of their Weyl symbols.

To end this review of preliminary properties of elliptic quadratic differential operators, let us underline that the property of normality in this class of operators can be easily
checked by computing the Poisson bracket of the real part and the imaginary part of their symbols
\begin{equation}\label{16}\inc
\{\textrm{Re }q, \textrm{Im } q\}=\frac{\partial \textrm{Re }q}{\partial \xi}.\frac{\partial \textrm{Im }q}{\partial x}-\frac{\partial \textrm{Re }q}{\partial x}.\frac{\partial \textrm{Im }q}{\partial \xi}. \num
\end{equation}

\bigskip

\begin{proposition}\label{17}
An elliptic quadratic differential operator 
$$q(x,\xi)^w : B \rightarrow \lde, \ n \in \nn^*,$$ 
is \emph{normal} if and only if the quadratic form defined by the Poisson bracket of the real part and the imaginary part of its symbol is equal to zero
\begin{equation}\label{17.1}\inc
\forall (x,\xi) \in \rr^{2n}, \ \{\emph{\textrm{Re }}q, \emph{\textrm{Im }} q\}(x,\xi)=0.\num
\end{equation}
\end{proposition}

\bigskip
\noindent
\textit{Proof of Proposition \ref{17}}. This proposition is a direct consequence of the composition formula in Weyl calculus (see Theorem 18.5.4 in \cite{hormander}), 
which induces that the Weyl symbol of the commutator 
$$[q^w, (q^w)^*]=[q^w, \overline{q}^w]=-2i [(\textrm{Re }q)^w, (\textrm{Im } q)^w],$$ 
is equal to 
$$-2i( \textrm{Re } q \ \sharp \ \textrm{Im } q - \textrm{Im } q \ \sharp \ \textrm{Re } q)=-2 \{\textrm{Re } q,\textrm{Im } q\},$$
because $\textrm{Re }q$ and $\textrm{Im }q$ are some quadratic forms. The notation $\textrm{Re } q \ \sharp \ \textrm{Im } q$ stands for the Weyl symbol of the operator 
obtained by composition $(\textrm{Re} q)^w(\textrm{Im} q)^w$. $\Box$ 
  
\bigskip

\noindent
\textit{Remark.} Let us notice that the symplectic invariance of the Poisson bracket (see (21.1.4) in \cite{hormander}), 
\begin{equation}\label{17.2}\inc
\{(\textrm{Re }q) \circ \chi, (\textrm{Im }q) \circ \chi\}=\{\textrm{Re }q,\textrm{Im }q\} \circ \chi, \num
\end{equation}
if $\chi$ stands for a linear symplectic transformation of $\rr^{2n}$, implies that the condition (\ref{17.1}) is symplectically invariant.

\subsection{Statement of the main results}
\init

Let us consider an elliptic quadratic differential operator 
$$q(x,\xi)^w : B \rightarrow L^2(\rr^n).$$
We know from (\ref{15}) that the spectrum of this operator is contained in the numerical range of its symbol $\Sigma(q)$. The following proposition gives a first localization
of the regions where the resolvent can blow up in norm and where spectral instabilities can occur.

\bigskip

\begin{proposition}\label{18}
Let $q : \rr^n \times \rr^n \rightarrow \cc$, $n \in \nn^*$, be a complex-valued elliptic quadratic form. We have 
$$\forall z \not\in \Sigma(q), \ \big\|\big{(}q(x,\xi)^w-z\big{)}^{-1}\big\| \leq \frac{1}{d\big{(}z,\Sigma(q)
\big{)}},$$
where $d\big{(}z,\Sigma(q)\big{)}$ stands for the distance from $z$ to the numerical range $\Sigma(q)$.
\end{proposition}

\bigskip
  
This result shows that the resolvent of an elliptic quadratic differential operator cannot blow up in norm far from the numerical range of its symbol. We are now going to 
study what kind of phenomena can occur in this particular set. There are two cases to separate according to the property of normality or non-normality of the operator.

\subsubsection{Case of a normal operator}

Let us consider a \textit{normal} elliptic quadratic differential operator 
$$q(x,\xi)^w : B \rightarrow L^2(\rr^n).$$
Let us recall that according to the proposition \ref{17} this property of normality is exactly equivalent to the fact that 
$$\forall (x,\xi) \in \rr^{2n}, \ \{\textrm{Re }q,\textrm{Im }q\}(x,\xi)=0.$$
In this case, we have the classical formula (\ref{1}) for its resolvent's norm
\begin{equation}\label{19}\inc
\forall z \not\in \sigma\big(q(x,\xi)^w\big), \ \big\|\big(q(x,\xi)^w-z\big)^{-1} \big\|=\frac{1}{d\big(z,\sigma(q(x,\xi)^w) \big)}, \num
\end{equation}
which induces that the $\eps$-pseudospectrum of this operator is exactly equal to the $\eps$-neighbourhood of its spectrum
$$\sigma_{\eps}\big(q(x,\xi)^w\big)=\big\{z \in \cc : d\big(z,\sigma( q(x,\xi)^w) \big) \leq \eps \big\}, \ \eps > 0.$$
This classical formula (\ref{19}) ensures that the resolvent cannot blow up in norm far from the spectrum and induces that the spectrum of such an operator is stable under small 
perturbations.

\bigskip

\noindent
\textit{Example 1.} The operator 
\inc
\begin{multline}\label{19.5}
q_1(x,\xi)^w=-(1+i)\partial_{x_1}^2-\partial_{x_2}^2+4(-1+i) x_1 \partial_{x_1}+2(-1+i) x_2 \partial_{x_1}+6i x_2 \partial_{x_2}
\\+2i x_1 \partial_{x_2}+(6+5i) x_1^2+(11+i)x_2^2+(10+4i)x_1 x_2-2+5i, \num
\end{multline}
is an example of a normal elliptic quadratic differential operator. Its spectrum is given by 
$$\sigma\big(q_1(x,\xi)^w \big)=\big\{(2k_1+1)+(2k_2+1) \sqrt{2}e^{i \frac{\pi}{4}} : (k_1,k_2) \in \nn^2 \big\}.$$
\begin{figure}[ht]
\caption{Spectrum and a $\eps$-pseudospectrum of the operator $q_1(x,\xi)^w$.}
\centerline{\includegraphics[scale=0.9]{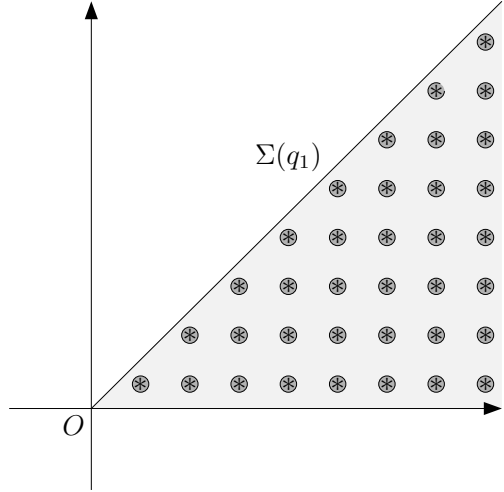}}
\end{figure}

\bigskip

\noindent
\textit{Example 2.} Let us notice that when the numerical range $\Sigma(q)$ is reduced to a closed half-line, the elliptic quadratic differential operator $q(x,\xi)^w$ is always
normal since
$$\{\textrm{Re }q,\textrm{Im }q\}=|z|^2\{\textrm{Re}(z^{-1}q),\textrm{Im}(z^{-1}q)\}=0,$$
if $z \in \cc^*$ is chosen such that $\textrm{Im}(z^{-1}q)=0$. In fact, the operator $q(x,\xi)^w$ can in this particular case be reduced after a conjugation by a unitary operator on 
$L^2(\rr^n)$ to the operator 
$$z \sum_{j=1}^{n}{\lambda_j (D_{x_j}^2+x_j^2)},$$
where $\lambda_j >0$ for all $j=1,...,n$.
\begin{figure}[ht]
\caption{Example of a normal elliptic quadratic differential operator.}
\centerline{\includegraphics[scale=0.9]{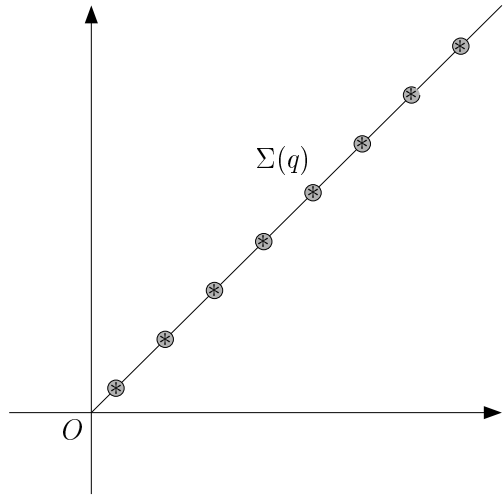}}
\end{figure}

\subsubsection{Case of a non-normal operator}\label{b2}
 
Let us consider a \textit{non-normal} elliptic quadratic differential operator
$$q(x,\xi)^w : B \rightarrow L^2(\rr^n), \ n \in \nn^*.$$
We assume in the following that the numerical range $\Sigma(q)$ is distinct from the whole complex plane
\begin{equation}\label{20}\inc
\Sigma(q) \neq \cc. \num
\end{equation}
As mentioned in the section \ref{b1}, this additional assumption is always fulfilled in dimension $n \geq 2$. 
It only excludes some very particular one-dimensional elliptic quadratic differential
operators (see the remark following the proposition \ref{29} for more precision about these operators).

Under this additional assumption, the numerical range $\Sigma(q)$ is always a closed angular sector with a top in $0$ and a \textit{positive} opening strictly lower than $\pi$.

\bigskip

\noindent
\textit{\ref{b2}.a. On the pseudospectrum at the interior of the numerical range.} Let us consider the associated semiclassical elliptic quadratic differential operator 
$$(q(x,h\xi)^w)_{0<h \leq 1}.$$ 
We can build in every point of the interior of the numerical range $\mathring{\Sigma}(q)$ some semiclassical quasimodes.

\bigskip

\begin{theorem}\label{21}
If the elliptic quadratic differential operator 
$$q(x,\xi)^w : B \rightarrow L^2(\rr^n), \ n \in \nn^*,$$ 
is \emph{non-normal} and verifies $\Sigma(q) \neq \cc$ then for all
$z \in \mathring{\Sigma}(q)$ and $N \in \nn$, there exist $h_0>0$ and a semiclassical family $(u_h)_{0<h \leq h_0} \in \mathcal{S}(\rr^n)$ such that 
$$\|u_h\|_{L^2(\rr^n)}=1 \textrm{ and } \|q(x,h\xi)^w u_h -z u_h\|_{L^2(\rr^n)}=O(h^N) \textrm{ when } h \rightarrow 0^+.$$
This result induces the existence of \emph{semiclassical pseudospectrum of infinite index} in every point of the interior of the numerical range $\mathring{\Sigma}(q)$.
\end{theorem}

\bigskip

According to (\ref{7}), this result in the semiclassical setting induces that the resolvent's norm of the quantum operator $q(x,\xi)^w$ blows up fastly along all the half-lines belonging 
to the interior of the numerical range $\mathring{\Sigma}(q)$,
\begin{equation}\label{22}\inc 
\forall z \in \mathring{\Sigma}(q), \forall N \in \nn, \forall C>0, \forall \eta_0 \geq 1, \exists \eta \geq \eta_0, \ \|\big(q(x,\xi)^w-z \eta \big)^{-1}\| \geq C \eta^N. \num
\end{equation}
We deduce from (\ref{15}) that as soon as an elliptic quadratic differential operator is \textit{non-normal} its resolvent blows up in norm in some regions of the resolvent set
far from its spectrum. This fact induces 
that the high energies of such an operator are very unstable under small perturbations as we have already noticed on the numerical computation performed for the 
rotated harmonic oscillator. It follows that in the class of elliptic quadratic differential operators\footnote{If we exclude the one-dimensional particular cases previously mentioned.} the 
property of spectral stability is exactly equivalent to the property of normality: 
$$\begin{array}{lllll}
\sigma( q(x,\xi)^w) \textrm{ is \textit{stable} under} & \Leftrightarrow & q(x,\xi)^w \textrm{ is a \textit{normal}} & \Leftrightarrow & \{\textrm{Re }q,\textrm{Im }q\}=0.\\
\textrm{small perturbations} & & \textrm{operator}& & 
\end{array}$$ 
By spectral stability, we mean here that the resolvent of these operators cannot blow up in norm far from their spectra. Let us add that it is not very surprising to have 
this property of spectral stability under the assumption of normality, but it is worth noticing that as soon as this property is violated, it occurs in this class of operators
some strong spectral instabilities under small perturbations for their high energies.

\bigskip

\noindent
\textit{Examples.} The two following operators 
\begin{equation}\label{23}\inc
q_2(x,\xi)^w=-\partial_{x_1}^2-2\partial_{x_2}^2+4i x_2 \partial_{x_2}+2x_1^2+(4+i) x_2^2+4x_1 x_2 +2i \num
\end{equation}
and\inc
\begin{multline}\label{24}
q_3(x,\xi)^w=-(1+i) \partial_{x_1}^2-2\partial_{x_2}^2+4(-1+i) x_1 \partial_{x_1}+2(1-i)x_2 \partial_{x_1}-4i x_1 \partial_{x_2}
\\+(9+4i)x_1^2+(2+i)x_2^2-4(1+i)x_1 x_2-2+2i, \num
\end{multline}
are some examples of non-normal elliptic quadratic differential operators.
 
\bigskip

\noindent
\textit{\ref{b2}.b. On the pseudospectrum at the boundary of the numerical range.} Let us now study what occurs on the boundary of the numerical range $\partial \Sigma(q)$ for 
a \textit{non-normal} elliptic quadratic differential operator
$$q(x,\xi)^w : B \rightarrow L^2(\rr^n).$$ 
Let us mention that we always assume that $\Sigma(q) \neq \cc$. Under these assumptions, the boundary of the numerical range is composed of the union of the origin $0$ and two half-lines
$\Delta_1$ and $\Delta_2$,
\begin{equation}\label{25}\inc
\partial \Sigma(q)=\{0\}  \sqcup \Delta_1 \sqcup \Delta_2, \num
\end{equation} 
that we can write 
\begin{equation}\label{26}\inc
\Delta_1=z_1 \rr_+^* \textrm{ and } \Delta_2=z_2 \rr_+^* \textrm{ with } z_1, z_2 \in \partial \Sigma(q) \setminus \{0\}. \num
\end{equation}
We need to define a notion of \textit{order} for the symbol $q(x,\xi)$ on these two half-lines $\Delta_j$, $j=1,2$. Let us begin by recalling the classical definition of the order $k(x_0,\xi_0)$
of a symbol $p(x,\xi)$ at a point $(x_0,\xi_0) \in \rr^{2n}$ (see section 27.2, chapter 27 in \cite{hormander}). This order $k(x_0,\xi_0)$ is an element of the set $\nn \cup \{+\infty\}$ defined by
\begin{equation}\label{t10.5}\inc
k(x_0,\xi_0)=\sup\big{\{}j \in \mathbb{Z} : p_I(x_0,\xi_0)=0, \ \forall\ 1 \leq |I| \leq j\big{\}}, \num
\end{equation}
where $I=(i_1,i_2,...,i_k) \in \{1,2\}^k$, $|I|=k$ and $p_I$ stands for the iterated Poisson brackets
$$p_{I}=H_{p_{i_1}}H_{p_{i_2}}...H_{p_{i_{k-1}}}p_{i_k},$$
where $p_1$ and $p_2$ are respectively the real and the imaginary part of the symbol $p$, $p=p_1+ip_2$. The order of a symbol $q$ at a point $z$ is then defined as the 
maximal order of the symbol $p=q-z$ at every point $(x_0,\xi_0) \in \rr^{2n}$ verifying 
$$p(x_0,\xi_0)=q(x_0,\xi_0)-z=0.$$ 
Let us underline that the symplectic invariance of the Poisson bracket (\ref{17.2}) induces the same property for the order of a symbol at a point.

Since here the symbol $q$ is a quadratic form, all the iterated Poisson brackets are also some quadratic forms. This property of degree two homogeneity of these Poisson brackets induces
that the symbol $q$ has the same order at every point of each half-line $\Delta_j$, $j=1,2$. This allows to define the order of the symbol $q$ on the half-line $\Delta_j$  by 
defining this order by this common value. Let us mention that this order can be \textit{finite} or \textit{infinite}.

\bigskip

\noindent
\textit{Examples.} One can easily check that the Weyl symbol
$$\xi^2+e^{i \theta}x^2, \ 0<\theta <\pi,$$
of the rotated harmonic oscillator has an order equal to $2$ on the both half-lines $\rr_+^*$ and $e^{i \theta}\rr_+^*$, which composes the boundary of its numerical range. The symbol
$q_2$ of the operator defined in (\ref{23}) has an order equal to $2$ on $i\rr_+^*$ and to $6$ on $\rr_+^*$,
$$\Sigma(q_2)=\{z \in \cc : \textrm{Re }z \geq 0, \ \textrm{Im }z \geq 0\}.$$
On the other hand, we can verify that the symbol $q_3$ of the operator defined in (\ref{24}) is of infinite order on the half-line $\rr_+^*$ and has an order equal to $2$ on $e^{i \pi/4}\rr_+^*$,
$$\Sigma(q_3)=\{0\} \cup \{z \in \cc^* : 0 \leq \textrm{arg } z \leq \pi/4\}.$$

\bigskip

In the case where the symbol is of \textit{finite} order on a half-line $\Delta_j$, $j=1,2$, we have the following result.

\bigskip

\begin{theorem}\label{27}
If the Weyl symbol $q(x,\xi)$ of a non-normal elliptic quadratic differential operator is of \emph{finite} order $k_j$ on the half-line 
$$\Delta_j, \ j \in \{1,2\}, \ \Delta_j \subset \partial \Sigma(q) \setminus \{0\},$$ 
then this order is necessary \emph{even} and there is \textbf{\emph{no}} semiclassical pseudospectrum of index $k_j/(k_j+1)$ on 
$\Delta_j$ for the associated semiclassical operator
$$\Delta_j \subset \cc \setminus \Lambda_{k_j/(k_j+1)}^{\emph{\textrm{sc}}}\big( q(x,h\xi)^w\big).$$
\end{theorem}

\bigskip

\noindent 
\textit{Remark.} Let us mention that we can more precisely establish that in dimension $n \geq 1$, the order $k_j$ is an even 
integer verifying
$$2 \leq k_j \leq 4n-2.$$
This result is proved in \cite{karel6}.

\bigskip

By rephrasing this result in a quantum setting, it follows from (\ref{8}) and (\ref{15}) that when the symbol $q$ of a non-normal elliptic quadratic differential operator $q(x,\xi)^w$
is of \textit{finite} order $k_j$ on a half-line 
$$\Delta_j, \ j \in \{1,2\}, \ \Delta_j \subset \partial \Sigma(q) \setminus \{0\},$$
then the resolvent of this operator remains bounded in norm in a set of the following type
\begin{equation}\label{28}\inc
\big{\{}u \in \cc :  |u| \geq C_1, \ d(\Delta_j,u) \leq C_2 |\textrm{proj}_{\Delta_j} u|^{\frac{1}{k_j+1}} \big{\}}, \num
\end{equation}
where $C_1$ and $C_2$ are some positive constants.

As we will see in its proof, this absence of semiclassical pseudospectrum is linked to some properties of subellipticity. Let us just underline for the moment that the index 
$k_j/(k_j+1)$, which appears in this result is exactly equal to the loss appearing in the subelliptic estimate hidden behind this result.

About the case of \textit{infinite} order, the situation is much more complicated. Nevertheless, we can first notice in this case that we cannot expect to prove a stronger result than an 
absence of semiclassical 
pseudospectrum of index 1. Indeed, we can easily check on the example of the operator $q_3(x,\xi)^w$ defined in (\ref{24}) that its spectrum is given by    
$$\sigma\big(q_3(x,\xi)^w \big)= \big\{(2 k_1+1) \sqrt{2}+(2 k_2+1)3^{\frac{1}{2}} 2^{\frac{1}{4}}e^{i \frac{ \pi}{8}}: (k_1,k_2) \in \nn^2  \big\}.$$
We recall that the spectrum of this operator is only composed of eigenvalues and that its symbol is of infinite order on $\rr_+^*$. It follows from the structure of the spectrum and 
(\ref{8}) that if there is no semiclassical pseudospectrum of infinite index in a point of the half-line $\rr_+^*$, there is necessary no semiclassical pseudospectrum 
of index $\mu$ with an index $\mu \geq 1$. In fact, we can prove by using a result of exponential decay in time for the norm of contraction semigroups generated 
by elliptic quadratic differential operators (see \cite{karel6}) that there is never some semiclassical pseudospectrum of index 1 on all these half-lines of infinite order. Let us mention that this 
result of exponential decay will not be proved here but it will be explained in the following how it induces the absence of semiclassical pseudospectrum of index 1.

\subsubsection{About the geometry of $\eps$-pseudospectra for elliptic quadratic differential operators}

Let us now explain what are the consequences of these results on the geometry of $\eps$-pseudospectra for elliptic quadratic differential operators. Let us begin by considering 
the one-dimensional case which is a bit particular. In dimension $n=1$, an elliptic quadratic differential operator can be reduced after a similitude and a conjugation by a unitary
operator to the harmonic oscillator or to the rotated harmonic oscillator.

\bigskip

\begin{proposition}\label{29}
Let us consider $q : \rr \times \rr \rightarrow \cc$ a complex-valued elliptic quadratic form such that $\Sigma(q) \neq \cc$. For all $h>0$, there exist a unitary operator \emph{(}more precisely
a metaplectic operator\emph{)} $U_h$ on $L^2(\rr)$, which is an automorphism of the spaces $\mathcal{S}(\rr)$ and $B$, $z \in \cc^*$ and $\theta \in [0,\pi[$ such that 
$$\forall h>0, \ q(x,h\xi)^w=z U_h\big((hD_x)^2+e^{i\theta}x^2 \big)U_h^{-1}.$$ 
\end{proposition}

\bigskip

\noindent
\textit{Remark.} In the case where $\Sigma(q)=\cc$, an elliptic quadratic differential operator $q(x,\xi)^w$ can be reduced after a similitude and a conjugation by a 
unitary operator on $L^2(\rr^n)$ to the operator defined in the Weyl quantization by the symbol
$$(\xi+i x)(\xi +\eta x) \textrm{ with } \eta \in \cc, \ \textrm{Im } \eta>0,$$
or
$$(\xi-i x)(\xi +\eta x) \textrm{ with } \eta \in \cc, \ \textrm{Im } \eta<0,$$
depending on the value of its Fredholm index, which is equal to $-2$ in the first case and to $2$ in the second one. 

\bigskip

As we will see in the following, this proposition allows us to reduce the study of a one-dimensional non-normal elliptic quadratic differential operator verifying 
$$\Sigma(q) \neq \cc,$$ 
to the one of the 
rotated harmonic oscillator 
$$H_{\theta}=D_x^2+e^{i\theta} x^2, \ 0<\theta <\pi.$$
Let us mention that the previous results (Theorem \ref{21} and Theorem \ref{27}) were already known in the particular case of the rotated harmonic oscillator. Indeed, the existence of 
semiclassical quasimodes inducing the presence of semiclassical pseudospectrum of infinite index in every point of the 
interior of the numerical range for the associated semiclassical operator, is a direct consequence of a result proved by E.B.~Davies in \cite{daviessemi} (Theorem 1). 
About the absence of semiclassical pseudospectrum of index $2/3$ on the boundary of the numerical range, this result has been proved for the rotated harmonic oscillator in 
\cite{karel3}\footnote{Let us recall that the value of the order is equal to 2 in this case.}.

As proved in \cite{karel3}, this absence of semiclassical pseudospectrum allows to give a proof of a conjecture stated by L.S. Boulton in \cite{boulton}. It deals with the 
geometry of $\eps$-pseudospectra for the rotated harmonic oscillator.
Let us now recall some facts about this conjecture and some results proved by L.S. Boulton in \cite{boulton}.

L.S. Boulton has first proved (Theorem 3.3 in \cite{boulton}) that the resolvent of the rotated harmonic oscillator blows up in norm along all a family of curves of the following form 
$$\eta \mapsto b \eta+e^{i\theta}\eta^p,$$
where $b$ and $p$ are some positive constants verifying $1/3<p<3$,  
\begin{equation}\label{30}\inc
\big{\|}\big{(}H_{\theta}-(b \eta + e^{i\theta}\eta^p)\big{)}^{-1}\big{\|} \rightarrow +\infty \ \textrm{when} \ \eta \rightarrow +\infty. \num
\end{equation}
On the other hand, he also proved that the resolvent of this operator remains bounded in norm on two half-stripes parallel to the half-lines $\rr_+$ or $e^{i\theta}\rr_+$. 
More precisely, he proved that there exist some positive constants $d$ and $M_d$ such that
\begin{equation}\label{31}\inc
\sup_{\eta \in \rr_+^*, \ 0 \leq b \leq d}{\big{\|}\big{(}H_{\theta}-(\eta+ib) \big{)}^{-1}\big{\|}} \leq M_d,  \num
\end{equation}
\begin{equation}\label{32}\inc
\sup_{\eta \in \rr_+^*, \ 0 \leq b \leq d}{\big{\|}\big{(}H_{\theta}-e^{i\theta}(\eta-ib) \big{)}^{-1}\big{\|}} \leq M_d. \num
\end{equation}
These bounds provide some information about the shape of $\varepsilon$-pseudospectra of the operator $H_{\theta}$. 
Indeed, L.S. Boulton has proved using these results that for all sufficiently small value of the positive parameter $\varepsilon$, the $\varepsilon$-pseudospectra of the rotated harmonic oscillator
is contained in the shaded set appearing on the following figure. The eigenvalues appear on this figure marked by some $\diamond$.
\begin{figure}[hhh]\label{33}
\caption{A first localization of the $\eps$-pseudospectra of the rotated harmonic oscillator.}
\centerline{\includegraphics[scale=0.9]{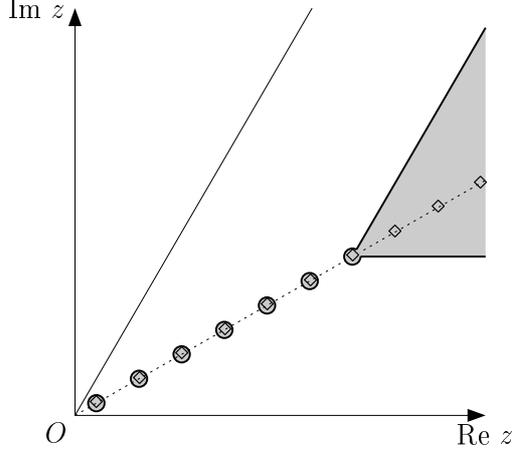}}
\end{figure}

More precisely, L.S. Boulton proved that for all
$0<\delta<1$ and $m \in
\mathbb{N}$, there exists a positive constant $\varepsilon_0$ such that for all $0<\varepsilon<\varepsilon_0$,
\begin{equation}\label{32}\inc
\sigma_{\varepsilon}(H_{\theta}) \subset \bigcup_{n=0}^{m}{\{z \in \cc : |z-\lambda_n|<\delta\}} \cup \big{[}\lambda_{m+1}-\delta
e^{i \theta/2}+S_{\theta}\big{]}, \num
\end{equation}
where 
$$\lambda_n=e^{i\theta/2}(2n+1), \ n \in \nn$$
and
$$S_{\theta}=\{z \in \cc^*:  0 \leq \textrm{arg}\ z \leq \theta\} \ \cup \ \{0\}.$$

\noindent
In fact, in view of some numerical calculations performed by E.B. Davies in \cite{daviesosc}, L.S.~Boulton has conjectured that the index $p=1/3$ appearing in (\ref{30}) is the
\textit{critical} one in the following sense:

Let us consider $0<p<1/3$, $0<\delta<1$ and $m \in \mathbb{N}$. If $b_{m,p}$ and~$E$ are some positive constants
verifying
$$b_{m,p} E+ e^{i\theta}E^p=\lambda_m \ \textrm{and} \ \forall \eta >E, \  \textrm{arg}\ z_{\eta}< \theta/2,$$
where $z_{\eta}=b_{m,p} \eta+e^{i\theta} \eta^p$, let us set
$$\Omega_{m,p}=\big{\{}|z_{\eta}|e^{i\alpha} \in \mathbb{C}: \ \eta \geq E, \ \textrm{arg}\ z_{\eta} \leq \alpha \leq 
\textrm{arg}(\overline{z_{\eta}}e^{i\theta})\big{\}}.$$
L.S. Boulton has conjectured the following result.

\bigskip

\noindent
\textbf{Boulton's conjecture}. There exists $\varepsilon_0>0$ such that for all $0<\varepsilon<\varepsilon_0$,
\begin{equation}\label{34}\inc
\sigma_{\varepsilon}(H_{\theta}) \subset \bigcup_{n=0}^{m}{\{z \in \mathbb{C}: \ |z-\lambda_n| <\delta\}} \cup \Omega_{m,p}.\num
\end{equation}

The absence of semiclassical pseudospectrum of index $2/3$ on the boundary of the numerical range $\partial \Sigma(q) \setminus \{0\}$ for the rotated harmonic 
oscillator\footnote{The order of the rotated harmonic oscillator's symbol is equal to 2 on $\partial \Sigma(q) \setminus \{0\}$.} given by the theorem \ref{27} shows that 
this index $1/3$ is actually the \textit{critical} one. Indeed, we can deduce (\ref{34}) from (\ref{28}) (see \cite{karel3} for more details) since here $k_j=2$, $j \in \{1,2\}$.
As we will see, this theorem \ref{27} is a consequence of a subelliptic estimate for general semiclassical pseudodifferential operators proved by 
N. Dencker, J. Sjöstrand and M.~Zworski in \cite{dencker} (Theorem 1.4). In the particular case of the rotated harmonic oscillator, a more elementary proof of this 
result using only some non-trivial localization scheme in the frequency variable is given in \cite{karel3}.

Let us notice that this inclusion (\ref{34}) allows to give a sharp description of the $\eps$-pseudospectra of the rotated harmonic oscillator, which is \textit{optimal} in view of (\ref{30}).
\begin{figure}[ht]
\caption{Shape of the $\eps$-pseudospectra of the rotated harmonic oscillator.}
\centerline{\includegraphics[scale=0.9]{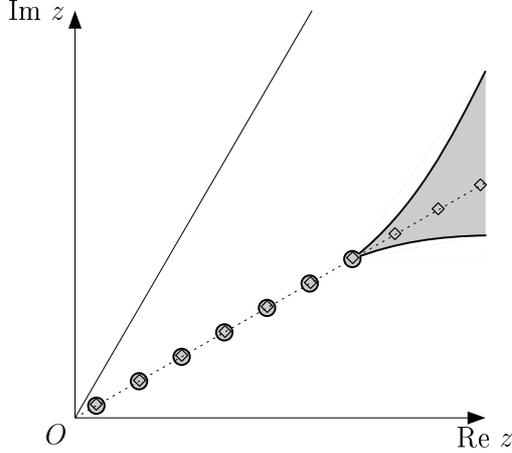}}
\end{figure}

By coming back to the case of an arbitrary dimension $n \geq 1$, let us finally underline that using the theorem \ref{27}, we can give similar descriptions of the $\eps$-pseudospectra 
for non-normal elliptic quadratic differential operators, to the one given by L.S. Boulton for the rotated harmonic oscillator, when the symbols of these operators are of \textit{finite} order 
on the two open half-lines, which compose the boundary of their numerical ranges. The only difference with the particular case of the rotated harmonic oscillator is that the critical indices, 
which appear in this description can be \textit{different}. Indeed, these critical indices depend directly according to (\ref{28}) on the order of the symbols on the two half-lines composing 
the boundary of their
numerical ranges. We refer the reader to \cite{karel3} for more details about the way of getting from (\ref{28}) such descriptions of $\eps$-pseudospectra.

\section{The proofs of the results}  
\init

Before giving the proofs of the results stated in the previous section, let us begin by recalling the \textit{symplectic invariance} property of the Weyl quantization (see Theorem 18.5.9 in 
\cite{hormander}). This symplectic invariance is actually the most important property of the Weyl quantization.

For every affine symplectic transformation $\chi$ of $\rr^{2n}$, there exists a unitary transformation $U$ on $\lde$, uniquely determined apart from a constant factor of 
modulus 1, such that $U$ is an automorphism of the spaces $\mathcal{S}(\rr^n)$, $B$ 
and $\mathcal{S}'(\rr^n)$, where $B$ is the Hilbert space defined in (\ref{14.1}), and
\begin{equation}\label{35}\inc
(a \circ \chi)(x,\xi)^w=U^{-1} a(x,\xi)^w U, \num
\end{equation}
for all $a \in \mathcal{S}'(\rr^{2n})$. The operator $U$ is a metaplectic operator associated to the affine symplectic transformation $\chi$.

This symplectic invariance of the Weyl quantization induces the same property for the semiclassical pseudospectra of elliptic quadratic differential operators in the sense that if 
$$q : \rr_x^n \times \rr_{\xi}^n \rightarrow \cc,$$ 
is a complex-valued elliptic quadratic form and $\chi$ is a \textit{linear} symplectic transformation of $\rr^{2n}$, we have for all $\mu \in [0,\infty]$,
\begin{equation}\label{36}\inc
\Lambda_{\mu}^{\textrm{sc}}\big( (q \circ \chi)(x,h\xi)^w\big)=\Lambda_{\mu}^{\textrm{sc}}\big( q(x,h\xi)^w\big). \num 
\end{equation}
To prove this fact, let us begin by noticing that for all $a \in \mathcal{S}'(\rr^{2n})$ and $h>0$, we have
$$U_h^{-1}a(x,\xi)^w U_h=a(h^{-1/2}x,h^{1/2}\xi)^w,$$
where 
$$U_h f(x)=h^{n/4}f(h^{1/2}x),$$
since according to the proof of Theorem 18.5.9 in \cite{hormander}, $U_h$ is a metaplectic operator 
associated to the linear symplectic transformation 
$$(x,\xi) \mapsto (h^{-1/2}x,h^{1/2}\xi).$$
Let us now consider the case where the symbol $a$ is a quadratic form. The homogeneity property of such a symbol implies that  
$$\forall h >0, \ a(h^{-1/2}x,h^{1/2}\xi)=\frac{1}{h}a(x,h\xi),$$
and  
$$\forall h>0, \ U_h^{-1} a(x,\xi)^w U_h=\frac{1}{h}a(x,h\xi)^w.$$
If $q : \rr_x^n \times \rr_{\xi}^n \rightarrow \cc$ is a complex-valued elliptic quadratic form and $\chi$ is a linear symplectic transformation of $\rr^{2n}$, we can notice that  
$$(q \circ \chi)(x,h\xi)^w, \ h>0,$$ 
is actually an elliptic quadratic differential operator since the symbol $q \circ \chi$ is an elliptic quadratic form. 
Let $z \in \cc$ and $U$ be a metaplectic operator associated to the linear symplectic transformation $\chi$. 
Using that $U$ and $U_h$ are some automorphisms of the Hilbert space $B$ and
\inc
\begin{multline*}\label{36.5}
U_h^{-1}U^{-1}U_h q(x,h\xi)^w U_h^{-1} U U_h=U_h^{-1}U^{-1} h q(x,\xi)^w U U_h\\
=h U_h^{-1} (q \circ \chi)(x,\xi)^w U_h=(q \circ \chi)(x,h\xi)^w, \num
\end{multline*} 
we obtain that 
$$U_h^{-1}U^{-1}U_h \big(q(x,h\xi)^w-z\big)^{-1} U_h^{-1} U U_h=\big((q \circ \chi)(x,h\xi)^w-z\big)^{-1}.$$
Using finally that $U_h^{-1} U^{-1} U_h$ is a unitary transformation of $L^2(\rr^n)$, this identity implies that    
$$\big\| \big((q \circ \chi)(x,h\xi)^w-z \big)^{-1}\big\|=\big\| \big( q(x,h\xi)^w-z\big)^{-1}\big\|,$$
which proves (\ref{36}).
In the following, this property of symplectic invariance will allow us to reduce certain symbols to some normal forms by choosing new symplectic coordinates. 
We can now begin to prove the results stated in the previous section. 

\bigskip

Let us start by the proof of the proposition \ref{18}.

\bigskip
\noindent
\textit{Proof of Proposition \ref{18}}.  
If the numerical range is equal to the whole complex plane, there is nothing to prove. If $\Sigma(q) \neq \cc$, we have seen in the previous section that 
the numerical range is necessary a closed angular sector with a top in $0$ and an opening strictly lower than $\pi$.

Let us consider $z \not\in \Sigma(q)$
and denote by $z_0$ its orthogonal projection on the non-empty closed convex set $\Sigma(q)$. According to the shape of the numerical range, it follows that $z_0$ 
belongs to its boundary and that we can find a complex number $z_1 \in \cc^*$, $|z_1|=1$ such that  
$$\Sigma(z_1 q) \subset \big\{z \in \cc : \textrm{Re } z \geq 0\big\}$$
and
\begin{equation}\label{37}\inc
z_1 z \in \big\{z \in \cc : \textrm{Re } z<0\big\}, \ d\big{(}z,\Sigma(q)\big{)}=d(z_1 z,i\rr). \num
\end{equation}
Using now that the operator $i [\textrm{Im}(z_1 q)]^w$ 
is formally skew-selfadjoint, we obtain that for all $u \in \mathcal{S}(\rr^n)$,
\inc
\begin{align*}\label{38}
& \ \textrm{Re}\big{(}z_1 q(x,\xi)^w u - z_1 z u,u \big{)}_{\lde}\\
=& \ d(z_1
z,i\rr)\|u\|_{\lde}^2+\big{(}\big[\textrm{Re}\big{(}z_1q(x,\xi)\big{)}\big]^w u,u\big{)}_{\lde}. \num
\end{align*}
Then, since the quadratic form $\textrm{Re}(z_1 q)$ is non-negative, we deduce from the symplectic invariance of the Weyl quantization and the theorem 21.5.3 in \cite{hormander} that there exists a metaplectic operator $U$ such that  
$$\big[\textrm{Re}\big(z_1 q(x,\xi)\big)\big]^w=U^{-1} \Big(\sum_{j=1}^{k}{\lambda_j(D_{x_j}^2+x_j^2)}
+\sum_{j=k+1}^{k+l}{x_j^2} \Big) U,$$
with $k,l \in \nn$ and $\lambda_j>0$ for all $j=1,...,k$. By using that $U$ is a unitary operator on $L^2(\rr^{n})$, we obtain that the quantity
\begin{align*}
& \ \big{(} \big[\textrm{Re}\big{(}z_1q(x,\xi)\big{)}\big]^w u,u\big{)}_{L^2(\rr^n)}\\
=& \ \sum_{j=1}^{k}{\lambda_j \big{(}\|D_{x_j} U u\|_{\lde}^2+\|x_j U u\|_{\lde}^2 \big{)}}+\sum_{j=k+1}^{k+l}{\|x_j U u\|_{\lde}^2},
\end{align*}
is non-negative. Then, we can deduce from the Cauchy-Schwarz inequality, $(\ref{37})$ and $(\ref{38})$ that for all $u \in \mathcal{S}(\rr^n)$,  
$$d\big{(}z,\Sigma(q) \big{)} \|u\|_{\lde} \leq |z_1| \ \|q(x,\xi)^w u -z u\|_{\lde}.$$
Finally, using the density of the Schwartz space $\mathcal{S}(\rr^n)$ in $B$ and the fact that $|z_1|=1$, we obtain that 
$$\forall z \not\in \Sigma(q), \ \big\|\big{(}q(x,\xi)^w-z\big{)}^{-1}\big\| \leq \frac{1}{d\big{(}z,\Sigma(q)\big{)}},$$
since according to (\ref{15}),
$\sigma\big( q(x,\xi)^w\big) \subset \Sigma(q). \ \Box$ 

\bigskip

We now consider the one-dimensional case, which is a bit particular.

\subsection{The one-dimensional case} 
\init

In dimension $n=1$, we can reduce the study of complex-valued elliptic quadratic forms to exactly three normal forms after a similitude and a real linear symplectic transformation.  

\bigskip

\begin{lemma}\label{39} 
Let $q : \rr_x \times \rr_{\xi} \rightarrow \cc$ be a complex-valued elliptic quadratic form in dimension $1$. Then, there exists a linear symplectic transformation 
$\chi$ of $\rr^2$ such that the symbol $q \circ \chi$ is equal to one of the following normal forms$:$
\medskip

\noindent
\emph{(}i\emph{)} \ \  $\alpha(\xi^2+e^{i \theta}x^2)$ with $\alpha \in \cc^*$, $0 \leq \theta < \pi$.\\
\emph{(}ii\emph{)} \ $\alpha(\xi+i x)(\xi+\eta x)$ with $\alpha \in \cc^*$, $\eta \in \cc$, $\emph{\textrm{Im }} \eta >0$.\\
\emph{(}iii\emph{)} $\alpha(\xi-i x)(\xi+\eta x)$ with $\alpha \in \cc^*$, $\eta \in \cc$, $\emph{\textrm{Im }} \eta <0$.

\medskip

\noindent
In the two last cases $(ii)$ and $(iii)$, the numerical range $\Sigma(q)$ is equal to the whole complex plane, $\Sigma(q)=\cc$. 
\end{lemma}

\bigskip
\noindent
\textit{Proof of Lemma \ref{39}.} Let $q :  \rr^2 \rightarrow \cc$ be a complex-valued elliptic quadratic form. Let us first consider the case where $\Sigma(q) \neq \cc$. We deduce from the
proposition~\ref{13} that we can reduce our study to the case where $\textrm{Re } q$ is a positive definite quadratic form. Then, using Lemma 18.6.4 in \cite{hormander},
we can find a real linear symplectic transformation to reduce the quadratic form $\textrm{Re } q$ to the normal form
$$\lambda(x^2+\xi^2), \textrm{ with } \lambda>0.$$
It follows that there exist some real constants $a,b$ and $c$ such that 
$$q(x,\xi)=\lambda \big{(}x^2+\xi^2+i (a x^2+ 2 b x \xi +c \xi^2) \big{)}.$$
Then, we can choose an orthogonal matrix $P \in O(2,\rr)$ diagonalizing the real symmetric matrix associated to the quadratic form $a x^2+ 2 b x \xi +c \xi^2$, 
$$P^{-1}\left(
  \begin{array}{cc}
  a & b \\
  b & c \\
  \end{array}
\right)P=\left(
  \begin{array}{cc}
  \lambda_1 & 0 \\
  0 & \lambda_2 \\
  \end{array}
\right),$$
with $\lambda_1, \lambda_2 \in \rr$. If $P \in O(2,\rr) \setminus SO(2,\rr)$, we have 
$$\tilde{P}^{-1}\left(
  \begin{array}{cc}
  a & b \\
  b & c \\
  \end{array}
\right)\tilde{P}=\left(
  \begin{array}{cc}
  \lambda_2 & 0 \\
  0 & \lambda_1 \\
  \end{array}
\right),$$
if $\sigma_0$ is the matrix with determinant equal to $-1$,
$$\left(
  \begin{array}{cc}
  0 & 1 \\
  1 & 0 \\
  \end{array}
\right),$$
and $\tilde{P}=P \sigma_0$. It follows that we can always diagonalize the real symmetric matrix associated to the quadratic form
$\lambda^{-1} \textrm{Im } q$ by conjugating it by an element of $SO(2,\rr)$. Since the symplectic group is equal in dimension 1 to the group 
$SL(2,\rr)$, we can after a linear symplectic transformation of $\rr^2$ reduce the quadratic form $q$ to  
$$\lambda\big{(}x^2+\xi^2+i(\gamma_1 x^2+\gamma_2 \xi^2) \big{)}=\alpha(\xi^2+ r e^{i \theta}x^2),$$ 
where $\gamma_1, \gamma_2 \in \rr$, $\alpha \in \cc^*$, $r>0$ and $\theta \in ]-\pi,\pi[$. Let us notice that the ellipticity of $q$ actually implies that 
$\theta \not\equiv \pi[2\pi]$. Finally, using the real linear symplectic transformation $(x,\xi) \mapsto (r^{-1/4} x, r^{1/4} \xi),$ we get a symbol of type $(i)$,
$$\alpha r^{1/2}(\xi^2+e^{i \theta} x^2),$$
if $0 \leq \theta < \pi$. If $-\pi < \theta <0$, we need to use besides the real linear symplectic transformation 
$(x,\xi) \mapsto (\xi,-x)$ to obtain a symbol of type $(i)$,
$$\alpha r^{\frac{1}{2}}e^{i \theta}(\xi^2+e^{-i \theta}x^2).$$ 
Let us now assume that $\Sigma(q)=\cc$. Since the dimension is equal to 1, we can factor the symbol $q$ on $\cc$ as a polynomial function of degree 2 in the variable~$\xi$. 
Thus, according to the dependence in the variable $x$ of the polynomial function's coefficients, we can find some complex numbers $\lambda_1, \lambda_2$ and $\alpha \in \cc^*$ such that 
$$q(x,\xi)=\alpha(\xi- \lambda_1 x)(\xi-\lambda_2 x).$$
The ellipticity assumption for the quadratic form $q$ induces that 
$$\textrm{Im } \lambda_j \neq 0,$$
if $j=1,2$. Using now the linear symplectic transformation 
$(x,\xi) \mapsto (x,\xi+\textrm{Re } \lambda_1 x),$
we can assume that
\begin{equation}\label{40}\inc
q(x,\xi)=\alpha(\xi-i r x)(\xi+ b x), \num
\end{equation}
with $r \in \rr^*$ and $\textrm{Im } b \neq 0$. Let us now check that the assumption $\Sigma(q)=\cc$ induces that $r \ \textrm{Im }b<0$. Since
$$(\xi-i r x)(\xi+b x)=\xi^2+(b-i r) x \xi-i r b x^2,$$
the condition $\Sigma(q)=\cc$ implies that for all $(v,w) \in \rr^2$, there exists a solution $(x_0,\xi_0) \in \rr^2$ of the system
\begin{equation}\label{41}\inc
\left\lbrace 
\begin{array}{l}
\xi^2+ \textrm{Re }b \ x \xi +r \ \textrm{Im }b \ x^2=v \\
x \xi (\textrm{Im }b-r)- r \ \textrm{Re }b \ x^2 =w. \num
\end{array}
\right.
\end{equation}
Let us first notice that the second equation of $(\ref{41})$ is fulfilled for all $w \in \rr$ only if 
$$\textrm{Im }b \neq r.$$ 
If $w \neq 0$, it follows from the 
second equation of $(\ref{41})$ that $x_0 \neq 0$ and
\begin{equation}\label{42}\inc
\xi_0=\frac{w+r \ \textrm{Re }b \ x_0^2}{(\textrm{Im }b-r)x_0}. \num
\end{equation} 
Let us consider the case where $v=0$. Using $(\ref{42})$ and the first equation of $(\ref{41})$, we obtain that 
$$(w+r \ \textrm{Re }b \ x_0^2)^2+ \textrm{Re }b \ (\textrm{Im }b-r) x_0^2(w+r \ \textrm{Re }b \ x_0^2)
+r \ \textrm{Im }b \ (\textrm{Im }b-r)^2 x_0^4=0.$$
We can rewrite this equation as $f_w(X_0)=0$ if we set $X_0=x_0^2$ and
\begin{equation}\label{43}\inc
f_w(X)=r \ \textrm{Im }b \ \big((\textrm{Re }b)^2+(\textrm{Im }b-r)^2 \big)X^2+w \ \textrm{Re }b \ (\textrm{Im }b+r)X+w^2. \num
\end{equation} 
Thus, the condition $\Sigma(q)=\cc$ implies that there exists for all $w \neq 0$, a non-negative solution $X_0$ of the equation $f_w(X_0)=0$.
Since the quantity $r \ \textrm{Im }b$ is assumed to be non-zero, we first study the case where $r \ \textrm{Im }b>0$. In this case, since 
\begin{equation}\label{44}\inc
f_w'(X)=2r \ \textrm{Im }b \ \big((\textrm{Re }b)^2+(\textrm{Im }b-r)^2 \big)X+w \ \textrm{Re }b \ (\textrm{Im }b+r) \num
\end{equation}   
and
$$2r \ \textrm{Im }b \ \big((\textrm{Re }b)^2+(\textrm{Im }b-r)^2 \big)>0,$$
because $\textrm{Im }b \neq r$, we have
\begin{equation}\label{45}\inc
\forall X \in \rr_+, \ f_w(X) \geq f_w(0)=w^2>0, \num
\end{equation}
if $w \neq 0$ and
$$-\frac{w \ \textrm{Re }b \ (\textrm{Im }b+r)}{2r \ \textrm{Im }b \ \big((\textrm{Re }b)^2+(\textrm{Im }b-r)^2 \big)} \leq 0.$$
The estimate $(\ref{45})$ shows that if $r \ \textrm{Im }b>0$, the equation $f_w(X)=0$ has no non-negative solution for all value of the parameter $w \neq 0$. 
This proves that the condition $\Sigma(q)=\cc$ induces that $r \ \textrm{Im }b<0$. 
Using the linear symplectic transformation 
$$(x,\xi) \mapsto (|r|^{-1/2} x,|r|^{1/2} \xi),$$
we obtain the normal forms  
$(ii)$ and $(iii)$,
$$\alpha |r|(\xi+i x)(\xi+ \eta x) \ \textrm{with} \ \textrm{Im } \eta>0 \ \textrm{and} \ \alpha |r|(\xi-ix)(\xi+\eta x) \ \textrm{with} \ \textrm{Im } \eta<0,$$
where $\eta=|r|^{-1}b.$ Finally, we can easily check that the numerical ranges of the normal forms $(ii)$ and $(iii)$ are actually equal to the whole complex plane $\cc$. $\Box$

\bigskip

Let us notice that the proposition \ref{29} and the remark following its statement are some direct consequences of the symplectic invariance property of the Weyl quantization 
(see (\ref{36.5})) and the previous lemma. We can add that as proved after the lemma~3.1 in~\cite{hypoelliptic}, the Fredholm indices of the one-dimensional elliptic quadratic 
differential operators with symbols 
of type $(i)$, $(ii)$ and $(iii)$ are respectively equal to $0$, $-2$ and~$2$.

As we have mentioned in the previous section, the results of Theorem \ref{21} and Theorem \ref{27} are already known in the particular case of the rotated harmonic oscillator. The existence of 
semiclassical quasimodes inducing the presence of semiclassical pseudospectrum of infinite index in every point of the 
interior of the numerical range for the associated semiclassical operator, is a direct consequence of a result proved by E.B.~Davies in~\cite{daviessemi} (Theorem 1) and; the absence of 
semiclassical pseudospectrum of index $2/3$ on the boundary of the numerical range has been proved for the rotated harmonic oscillator in 
\cite{karel3}\footnote{Let us recall that the value of the order is equal to 2 in this case.}. 
As we have previously mentioned (see (\ref{17.2}) and (\ref{36})), the property of non-normality, the order of symbols and the semiclassical pseudospectra of elliptic quadratic differential operators
are \textit{symplectically invariant}. These properties allow us to reduce by any real linear symplectic transformations the symbols of the elliptic quadratic differential operators that we consider 
in our proof of the theorem \ref{21} and the theorem \ref{27}. By using the lemma \ref{39}, we deduce from the results of the theorem \ref{21} and the theorem \ref{27} proved for the 
rotated harmonic oscillator that they are therefore also fulfilled by all non-normal one-dimensional elliptic quadratic differential operators with a numerical range different 
from the whole complex plane.

\bigskip

We now consider the multidimensional case. As we will see in the following, there is a real jump of complexity between the one-dimensional case and the multidimensional one.
This jump is among other things a consequence of the complexity increase of symplectic geometry in dimension $n \geq  2$ and the larger diversity appearing in the class 
of elliptic quadratic differential operators.

\subsection{Case of dimension~\mathversion{bold}$n \geq 2$}
\init

We only need to study the case of a \textit{non-normal} elliptic quadratic differential operator
\begin{equation}\label{45.8}\inc
q(x,\xi)^w : B \rightarrow L^2(\rr^n), \num
\end{equation}
in dimension $n \geq 2$. Let us recall that in this case, the numerical range $\Sigma(q)$ is a closed angular sector with a top in $0$ and a \textit{positive} opening 
strictly lower than $\pi$, and that the proposition \ref{17} gives that 
\begin{equation}\label{46}\inc
\exists (x_0,\xi_0) \in \rr^{2n}, \ \{\textrm{Re }q,\textrm{Im } q\}(x_0,\xi_0) \neq 0.\num
\end{equation}
Let us begin by studying what occurs at the interior of the numerical range $\mathring{\Sigma}(q)$.

\subsubsection{On the pseudospectrum at the interior of the numerical range}\label{46.5}

To prove the existence of semiclassical quasimodes for the associated semiclassical operator given by the theorem \ref{21}, we need a first purely algebraic step
to characterize the points belonging to the interior of the numerical range.

Let us consider the following decomposition of the numerical range
\begin{equation}\label{47}\inc
\Sigma(q)=\tilde{A} \sqcup \tilde{B}, \num
\end{equation}
where 
\begin{equation}\label{48}\inc
\tilde{A}=\big{\{}z \in \Sigma(q) : \exists (x_0,\xi_0) \in \rr^{2n}, \ z=q(x_0,\xi_0), \ \{\textrm{Re } q,\textrm{Im }q\}(x_0,\xi_0)\neq 0\big{\}} \num
\end{equation}
and
\begin{equation}\label{49}\inc
\tilde{B}=\big{\{}z \in \Sigma(q) : z=q(x_0,\xi_0) \Rightarrow \{\textrm{Re } q,\textrm{Im } q\}(x_0,\xi_0)=0\big{\}}.\num
\end{equation}
The next section is devoted to give a geometrical description of these two sets. We establish using purely algebraic arguments that  
\begin{equation}\label{50}\inc
\tilde{A}=\mathring{\Sigma}(q) \ \textrm{and } \tilde{B}=\partial\Sigma(q).\num
\end{equation}
This result is a consequence of the geometry induced by the quadratic setting to which the studied symbols belong.

Let us begin by noticing that the symplectic invariance of the Poisson bracket (\ref{17.2}) induces the same property for the sets $\tilde{A}$ and $\tilde{B}$. We can therefore 
use some real linear symplectic transformation to reduce the symbol $q$. Since 
$$\{\textrm{Re}(zq),\textrm{Im}(zq)\}=|z|^2\{\textrm{Re }q,\textrm{Im }q\},$$ 
we deduce from this symplectic invariance, from the proposition \ref{13} and the lemma 18.6.4 in \cite{hormander} that after a similitude, we can reduce our study to the case where
\begin{equation}\label{51}\inc
\textrm{Re }q(x,\xi)=\sum_{j=1}^{n}{\lambda_j(\xi_j^2+x_j^2)}, \num
\end{equation}
with $\lambda_j >0$ for all $j=1,...,n$.

\subsubsection*{\ref{46.5}.a. Geometrical description of the sets $\tilde{A}$ and $\tilde{B}$}

We begin by proving the following inclusion  
\begin{equation}\label{52}\inc
\partial \Sigma(q) \subset \tilde{B}. \num
\end{equation}

Let us consider $z \in \partial \Sigma(q)$ and $(x_0,\xi_0) \in \rr^{2n}$ such that $z=q(x_0,\xi_0)$. This is possible because the numerical range is a closed angular sector.  
If $z=0$, the ellipticity property of $q$ implies that 
$$(x_0,\xi_0)=(0,0) \textrm{ and } \{\textrm{Re } q,\textrm{Im } q\}(x_0,\xi_0)=0,$$
because this Poisson bracket is also a quadratic form. This proves that $z \in \tilde{B}$.
If 
$$z \in \partial\Sigma(q) \setminus \{0\},$$ 
let us consider the global solution $Y$ of the linear Cauchy problem
\begin{equation}\label{53}\inc
\left\lbrace \begin{array}{l}
Y'(t)=H_{\textrm{Re } q}\big{(}Y(t)\big{)} \\
Y(0)=(x_0,\xi_0),
\end{array}\right. \num
\end{equation}
associated to the Hamilton vector field of the symbol $\textrm{Re } q$,
$$H_{\textrm{Re } q}=\sum_{j=1}^{n}{\Big(\frac{\partial \textrm{Re } q}{\partial \xi_j}\frac{\partial}{\partial x_j}- \frac{\partial \textrm{Re } q}{\partial x_j}\frac{\partial}{\partial \xi_j}\Big)}.$$
It is actually a linear Cauchy problem since $\textrm{Re }q$ is a quadratic form. Setting  
$$f(t)=\textrm{Im } q\big{(}Y(t)\big{)},$$ 
a direct computation gives that 
$$f'(0)=\{\textrm{Re } q,\textrm{Im } q\}(x_0,\xi_0).$$
If $f'(0) \neq 0$, we could find $t_0 \neq 0$ such that  
$$|f(t_0)| > |f(0)|=|\textrm{Im } z|.$$
Since $Y$ is the flow associated to the Hamilton vector field of $\textrm{Re }q$, the quadratic form $\textrm{Re } q$ is constant under it.
It follows that for all $t \in \rr$, 
$$\textrm{Re }q\big{(}Y(t)\big{)}=\textrm{Re }q\big{(}Y(0)\big{)}=\textrm{Re } z$$
and provides a contradiction because, since $z \in \partial\Sigma(q) \setminus \{0\}$, this would imply in view of the shape of the numerical range $\Sigma(q)$ (see Figure 7)
that 
$$q\big{(}Y(t_0)\big{)} \not\in \Sigma(q).$$ 
\begin{figure}[ht]
\caption{}
\centerline{\includegraphics[scale=0.9]{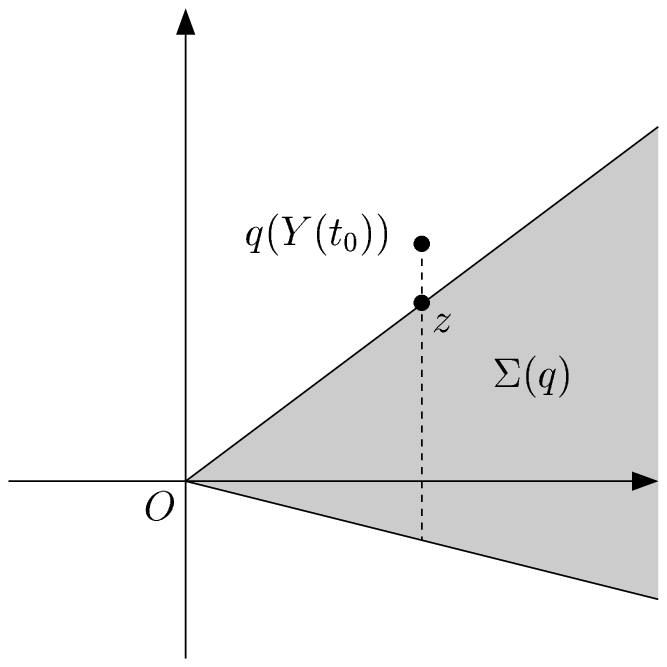}}
\end{figure}
\noindent
It follows that the Poisson bracket $\{\textrm{Re } q,\textrm{Im } q\}(x_0,\xi_0)$ is necessary equal to 0 and that $z \in \tilde{B}$. This ends the proof of the inclusion (\ref{52}).

\bigskip

Let us now assume that 
\begin{equation}\label{54}\inc
\partial \Sigma(q) \subset \tilde{B}, \ \partial \Sigma(q) \neq \tilde{B}. \num
\end{equation} 
In this case, we could find 
\begin{equation}\label{55}\inc
z \in \tilde{B} \setminus \partial \Sigma(q). \num
\end{equation}
Let us first notice that $z$ is necessary non-zero since $0 \in \partial\Sigma(q)$, and that $\textrm{Re } z>0$, since from (\ref{51}), 
\begin{equation}\label{55.1}\inc
\Sigma(q) \setminus \{0\} \subset \{z \in \cc^* : \textrm{Re } z>0\}. \num
\end{equation}
The fact that $z$ belongs to the set $\tilde{B}$ implies that 
\begin{equation}\label{56}\inc
\left\lbrace \begin{array}{l}
\textrm{Re } q(x,\xi)=\textrm{Re } z \\
\textrm{Im } q(x,\xi)=\textrm{Im } z
\end{array}\right. \Longrightarrow \{\textrm{Re } q,\textrm{Im } q\}(x,\xi)=0. \num
\end{equation}
We also know that there exists at least one solution to the system appearing in the left-hand-side of $(\ref{56})$. 
Since from (\ref{51}), the quadratic form $\textrm{Re }q$ is positive definite, we can simultaneously reduce the quadratic forms $\textrm{Re } q$ and $\textrm{Im } q$ 
by finding an isomorphism $P$ of $\rr^{2n}$ such that in the new coordinates $y=P^{-1}(x,\xi)$,
\begin{equation}\label{2.3.23}\inc
\textrm{Re } q(Py)=\sum_{j=1}^{2n}{y_j^2} \ \textrm{and} \ \textrm{Im } q(Py)=\sum_{j=1}^{2n}{\alpha_j y_j^2} \ \textrm{with} \ 
\alpha_1 \leq ... \leq \alpha_n. \num
\end{equation}
Let us now consider the following quadratic form  
\begin{equation}\label{2.3.24}\inc
p(y)=\{\textrm{Re } q,\textrm{Im } q\}(Py). \num
\end{equation}
We get from $(\ref{56})$ and $(\ref{2.3.23})$ that 
\begin{equation}\label{2.3.25}\inc
\left\lbrace \begin{array}{l}
\sum_{j=1}^{2n}{y_j^2}=\textrm{Re } z \\
\sum_{j=1}^{2n}{\alpha_j y_j^2}=\textrm{Im } z
\end{array}\right. \Longrightarrow p(y)=0. \num
\end{equation}
Let us underline that the isomorphism $P$ is not a priori a symplectic transformation and that it does not preserve the Poisson bracket $\{\textrm{Re } q,\textrm{Im } q\}$.

We consider the two following sets
\begin{equation}\label{2.3.26}\inc
E_1=\big\{y \in \rr^{2n} : r(y)=0\big\}, \num
\end{equation}
where
\begin{equation}\label{2.3.27}\inc
r(y)=\sum_{j=1}^{2n}{\Big(\alpha_j-\frac{\textrm{Im } z}{\textrm{Re } z} \Big) y_j^2} \num
\end{equation}
and
\begin{equation}\label{2.3.28}\inc
E_2=\big\{y \in \rr^{2n} : p(y)=0\big\}. \num
\end{equation}
The next lemma gives a first inclusion between these two sets $E_1$ and $E_2$.

\bigskip

\begin{lemma}\label{2.3.2}
We have 
\begin{equation}\label{2.3.29}\inc
E_1 \subset E_2. \num
\end{equation}
\end{lemma}

\bigskip

\noindent
\textit{Proof of Lemma \ref{2.3.2}}. Let $y \in E_1$. If $y=0$ then $y$ belongs to $E_2$ since from $(\ref{2.3.24})$, $p$ is a quadratic form in the variable $y$. 
If $y \neq 0$, we set 
$$t=\sum_{j=1}^{2n}{y_j^2}>0 \ \textrm{and} \ \forall j=1,...,2n, \ \tilde{y}_j=\sqrt{\frac{\textrm{Re } z}{t} }y_j.$$
We recall from $(\ref{55.1})$ that $z \in \tilde{B} \setminus \partial \Sigma(q)$ implies that $\textrm{Re } z>0$.
Then, since, on one hand 
$$\sum_{j=1}^{2n}{\tilde{y}_j^2}=\textrm{Re } z,$$
and that, on the other hand, we have from $(\ref{2.3.26})$ and $(\ref{2.3.27})$ that 
$$\sum_{j=1}^{2n}{\alpha_j \tilde{y}_j^2}=\frac{\textrm{Re } z}{t}\sum_{j=1}^{2n}{\alpha_j y_j^2}=
\frac{\textrm{Re } z}{t} \sum_{j=1}^{2n}{\frac{\textrm{Im } z}{\textrm{Re } z}y_j^2}= \textrm{Im } z,
$$
because $y \in E_1$, we deduce from $(\ref{2.3.25})$ and the homogeneity of degree 2 of the quadratic form $p$ that  
$$p(\tilde{y})=\frac{\textrm{Re } z}{t} p(y)=0.$$
According to $(\ref{2.3.28})$, this proves that $y \in E_2$ and ends the proof of the lemma~\ref{2.3.2}.~$\Box$

\bigskip

Then, we can notice from $(\ref{2.3.23})$ that the boundary of the numerical range $\partial \Sigma(q)$ is given by 
\begin{equation}\label{2.3.30}\inc
(1+i\alpha_1) \rr_+ \cup (1+i \alpha_n)\rr_+. \num
\end{equation}
Since the numerical range $\Sigma(q)$ is a closed set, the assumption  
$$z \in \tilde{B} \setminus \partial \Sigma(q) \subset \Sigma(q) \setminus \partial \Sigma(q) =\mathring{\Sigma}(q),$$ 
induces from $(\ref{2.3.30})$ that 
$$\frac{\textrm{Im } z}{\textrm{Re } z} \in ]\alpha_1,\alpha_n[.$$
This implies that the signature $(r_1,s_1)$ of the quadratic form $r$ defined in $(\ref{2.3.27})$ fulfills 
\begin{equation}\label{2.3.31}\inc
(r_1,s_1) \in \nn^* \times \nn^* \ \textrm{and} \ r_1+s_1 \leq 2n. \num
\end{equation}
Thus, we can assume after a new labeling that  
\begin{equation}\label{2.3.32}\inc
r(y)=a_1 y_1^2 +...+ a_{r_1} y_{r_1}^2-a_{r_{1}+1} y_{r_1+1}^2-...-a_{r_1+s_1} y_{r_1+s_1}^2, \num
\end{equation}
with $a_j>0$ for all $j=1,...,r_1+s_1$. It follows from $(\ref{2.3.26})$ and $(\ref{2.3.32})$ that in these new coordinates, the set $E_1$ is the direct product of a proper cone $C$ of $\rr^{r_1+s_1}$ 
and $\rr^{2n-r_1-s_1}$,
\begin{equation}\label{2.3.33}\inc
E_1=C \times \rr^{2n-r_1-s_1}. \num
\end{equation} 
\begin{figure}[ht]
\caption{}
\centerline{\includegraphics[scale=0.9]{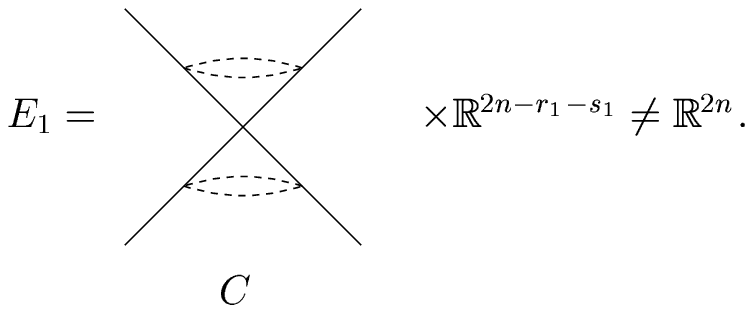}}
\end{figure}

We are now going to prove that the two sets $E_1$ and $E_2$ are equal
\begin{equation}\label{2.3.34}\inc
E_1=E_2. \num
\end{equation}
Let us reason by the absurd by assuming that it is not the case. Then, we could find from the lemma \ref{2.3.2},  
\begin{equation}\label{marise1}\inc
y_0 \in E_2 \setminus E_1, \ y_0=(y_0',y_0'') \textrm{ with } y_0' \in \rr^{r_1+s_1}, \  y_0'' \in \rr^{2n-r_1-s_1}. \num
\end{equation} 
We deduce from $(\ref{2.3.33})$ that $y_0' \not\in C$. Let us now recall an elementary geometrical fact that we will use several times. This fact is that the intersection of a real line 
and a real quadric surface is reduced to either $0,1$ or $2$ points, or the line is completely contained in the quadric surface. 
We first begin by proving that   
\begin{equation}\label{2.3.35}\inc
\rr^{r_1+s_1} \times \{y''=y_0''\} \subset E_2. \num
\end{equation}
Indeed, let us consider the affine subspace
$$F=\{y \in \rr^{2n} : y=(y',y'') \in \rr^{r_1+s_1} \times \rr^{2n-r_1-s_1}, \ y''=y_0''\}.$$ 
We identify for more simplicity the space $F$ to the space $\rr^{r_1+s_1}$. 
We agree to say that a point $x_0'$ of $\rr^{r_1+s_1}$ belongs to the set $E_2$ to mean that the point $(x_0',y_0'')$ belongs to the set $E_2$. 
With this convention, it is sufficient for proving the inclusion $(\ref{2.3.35})$ to consider some particular lines of $\rr^{r_1+s_1}$, containing the point $y_0'$ defined in (\ref{marise1}) and, which have an intersection 
with the cone $C$ in at least two other different points $u_0'$ and $v_0'$ (see Figure 9). These lines are necessary contained in the quadric surface $E_2$ because 
from the lemma \ref{2.3.2}, $$E_1 \subset E_2,$$ 
and that there are at least three different points of intersection between these lines and the quadric surface $E_2$,
$$(u_0',y_0'') \in C \times \rr^{2n-r_1-s_1}=E_1 \subset E_2, \ (v_0',y_0'') \in C \times \rr^{2n-r_1-s_1}=E_1 \subset E_2,$$
and $(y_0',y_0'') \in E_2$. 
Thus, we prove that the shaded disc appearing on the figure 10 is completely contained in the set $E_2$. 
By using the cone structure of the set $E_2$, we can deduce that all the interior of the cone $C$ (see Figure 11) is contained in $E_2$. 
Then, using again other particular intersections with some lines as on the figure 12, we deduce from our identification of the space $F$ to $\rr^{r_1+s_1}$ that 
the inclusion $(\ref{2.3.35})$ is fulfilled. 
\begin{figure}[ht]
\caption{}
\centerline{\includegraphics[scale=0.9]{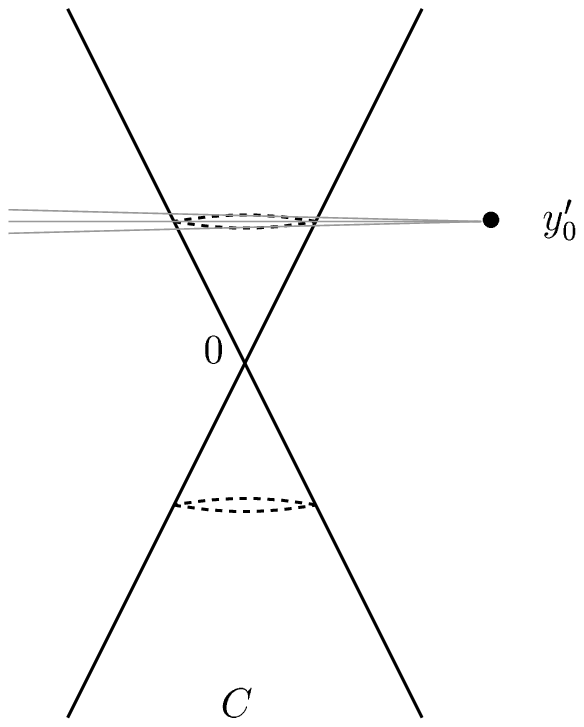}}
\end{figure}
\begin{figure}[ht]
\caption{}
\centerline{\includegraphics[scale=0.9]{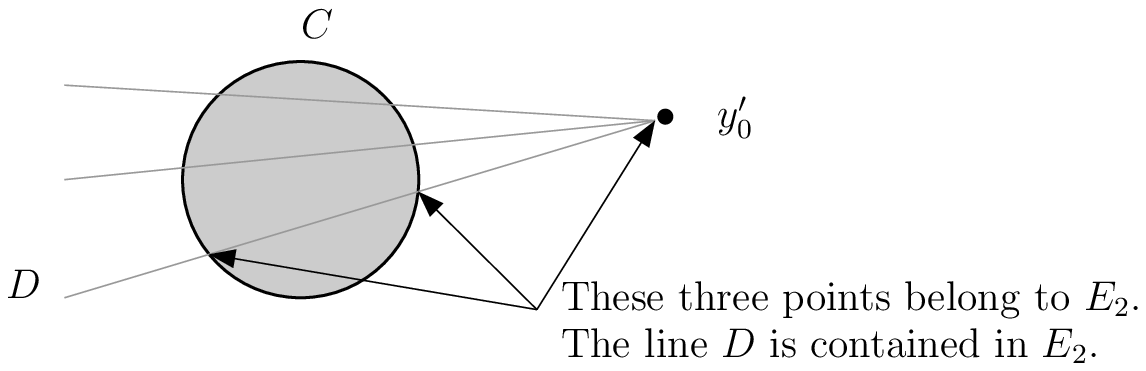}}
\end{figure}
\begin{figure}[ht]
\caption{}
\centerline{\includegraphics[scale=0.9]{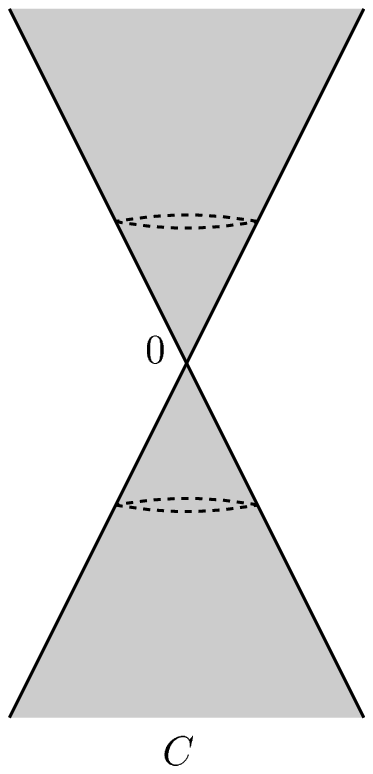}}
\end{figure}
\begin{figure}[ht]
\caption{}
\centerline{\includegraphics[scale=0.9]{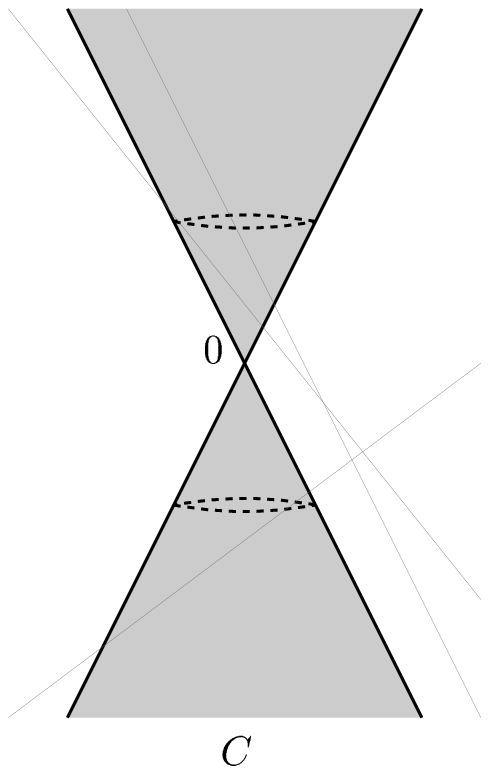}}
\end{figure}

We now prove that under these conditions, we have the identity 
\begin{equation}\label{2.3.36}\inc
E_2=\rr^{2n}. \num
\end{equation}
Indeed, let us consider $(\tilde{y}_0',\tilde{y}_0'') \in \rr^{2n}=\rr^{r_1+s_1}\times \rr^{2n-r_1-s_1}$. 
If $\tilde{y}_0' \in C$, then 
$$(\tilde{y}_0',\tilde{y}_0'') \in E_2,$$
because from $(\ref{2.3.29})$ and $(\ref{2.3.33})$, $(\tilde{y}_0',\tilde{y}_0'') \in E_1$ and $E_1 \subset E_2$. 
If, on the other hand $\tilde{y}_0' \not\in C$, we can choose a point $u \in \rr^{r_1+s_1}$ different from $\tilde{y}_0'$ such that $u \not\in C$, and such that the line containing $\tilde{y}_0'$ and 
$u$ in $\rr^{r_1+s_1}$, has an intersection with $C$ in at least two other different points $v$ and $w$ (see Figure 13).
\begin{figure}[ht]
\caption{}
\centerline{\includegraphics[scale=0.9]{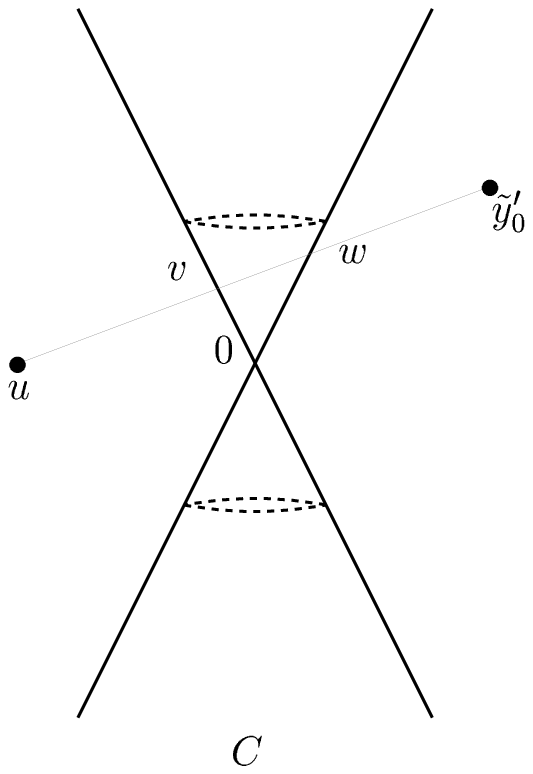}}
\end{figure}
Thus, we can find some distinct real numbers $t_1,t_2 \in \rr \setminus \{0,1\}$ such that 
$$v=(1-t_1) \tilde{y}_0'+ t_1 u \in C \ \textrm{and} \ w=(1-t_2) \tilde{y}_0'+ t_2 u \in C.$$
Considering now the line 
$$D=\big{\{}(1-t)(\tilde{y}_0',\tilde{y}_0'')+t (u,y_0'') : t \in \rr\big{\}},$$
we can notice that this real line contains at least three different points of $E_2$:
$$(v,(1-t_1)\tilde{y}_0''+t_1 y_0''), \ (w,(1-t_2)\tilde{y}_0''+t_2 y_0'') \ \textrm{and} \ (u,y_0'').$$
Indeed, this is a consequence of the fact that $v$ and $w$ belong to $C$, and from $(\ref{2.3.29})$, $(\ref{2.3.33})$ and $(\ref{2.3.35})$.
Thus, the line $D$ is contained in the quadric surface $E_2$. This implies that $(\tilde{y}_0',\tilde{y}_0'') \in D \subset E_2$.

To sum up, we have proved that if the two sets $E_1$ and $E_2$ are different then the set $E_2$ is equal to $\rr^{2n}$. This fact induces in view of $(\ref{2.3.28})$ that 
the quadratic form~$p$ is identically equal to zero. By coming back to the first coordinates $(x,\xi)=Py$, it follows from $(\ref{2.3.24})$ that the quadratic form 
$\{\textrm{Re } q, \textrm{Im } q\}$ is also identically equal to zero, which contradicts (\ref{46}).
This proves the identity $(\ref{2.3.34})$, 
$$E_1=E_2.$$
With this fact, we can resume our first reasoning by the absurd, which assume in $(\ref{55})$ the existence of a point $z \in \tilde{B} \setminus \partial \Sigma(q)$. 
Let us now consider $y_0 \not\in E_1=E_2$. This is possible according to $(\ref{46})$, $(\ref{2.3.24})$ and $(\ref{2.3.28})$. 
We deduce from $(\ref{2.3.26})$ and $(\ref{2.3.28})$ that $r(y_0)$ and $p(y_0)$ are non-zero. 
By considering $\lambda \in \rr^*$ such that 
$$p(y_0)=\lambda r(y_0)$$ 
and 
\begin{equation}\label{2.3.37}\inc
\tilde{r}(y)=p(y)-\lambda r(y), \num
\end{equation}
it follows from $(\ref{2.3.26})$, $(\ref{2.3.28})$, $(\ref{2.3.34})$ and $(\ref{2.3.37})$ that  
\begin{equation}\label{2.3.38}\inc
E_1 \subset \{y \in \rr^{2n} : \tilde{r}(y)=0\}. \num
\end{equation}
This inclusion $(\ref{2.3.38})$ is strict since 
$$\tilde{r}(y_0)=0 \textrm{ and } y_0 \not\in E_1.$$ 
By using now exactly the same reasoning as the one previously described to prove $(\ref{2.3.34})$, about the intersections of real lines and quadric surfaces,
we prove that the quadratic form $\tilde{r}$ is necessary identically equal to zero. Then, it follows from $(\ref{2.3.37})$ that 
\begin{equation}\label{2.3.39}\inc
p=\lambda r. \num
\end{equation}
By coming back to the first coordinates $(x,\xi)=Py$, we get using $(\ref{2.3.23})$, $(\ref{2.3.24})$, $(\ref{2.3.27})$ and $(\ref{2.3.39})$ that for all $(x,\xi) \in \rr^{2n}$, 
\begin{equation}\label{2.3.40}\inc
\{\textrm{Re } q, \textrm{Im } q\}(x,\xi)=\lambda \Big(\textrm{Im } q(x,\xi)- \frac{\textrm{Im }
z}{\textrm{Re } z}\ \textrm{Re } q(x,\xi) \Big). \num
\end{equation}
Let us now consider $(x_0,\xi_0) \in \rr^{2n}$ such that $q(x_0,\xi_0) \in \partial \Sigma(q) \setminus \{0\}$. This is possible since the numerical range $\Sigma(q)$ 
is a closed angular sector with a top in $0$ and a positive opening. We deduce from $(\ref{49})$ and $(\ref{52})$ that we necessarily have  
$$\{\textrm{Re } q,\textrm{Im } q\}(x_0,\xi_0)=0.$$ 
This induces from $(\ref{2.3.40})$ that 
\begin{equation}\label{2.3.41}\inc
\textrm{Im } q(x_0,\xi_0)=\frac{\textrm{Im } z}{\textrm{Re } z}\ \textrm{Re } q(x_0,\xi_0), \num
\end{equation}
because $\lambda \in \rr^*$. Since according to the shape of the numerical range $\Sigma(q)$ and $(\ref{55.1})$, 
$$q(x_0,\xi_0) \in \partial \Sigma(q) \setminus \{0\} \subset \{z \in \cc : \textrm{Re } z>0\},$$ 
the identity (\ref{2.3.41}) proves that the point $z$ also belongs to the set $\partial \Sigma(q)$, but it contradicts the initial assumption
$$z \in \tilde{B} \setminus \partial \Sigma(q).$$ 
Finally, this ends our reasoning by the absurd and proves (\ref{50}).

\subsubsection*{\ref{46.5}.b. Existence of semiclassical quasimodes at the interior of the numerical range}

To prove the existence of semiclassical quasimodes for the associated semiclassical operator 
$$(q(x,h\xi)^w)_{0<h \leq 1},$$
in every point of the numerical range's interior (Theorem \ref{21}), we use an existence result of semiclassical quasimodes for general pseudodifferential operators violating 
the condition $(\overline{\Psi})$\footnote{The definition of the condition $(\overline{\Psi})$ is recalled below.}. Let us mention that this result generalizes 
the two existence results of semiclassical quasimodes given by E.B. Davies, in the case of Schrödinger operators (Theorem 1 in \cite{daviessemi}), and by M.~Zworski 
in \cite{zworski1} and \cite{zworski2}, for pseudodifferential operators.


This existence result of semiclassical quasimodes can be stated as follows. 
Let us consider a semiclassical symbol $P(x,\xi;h)$ in $S(\langle (x,\xi) \rangle^m,dx^2+d\xi^2 )$ with $m \in \rr_+$, 
$$\langle (x,\xi) \rangle^2=1+x^2+\xi^2,$$ 
where $S(\langle (x,\xi) \rangle^m,dx^2+d\xi^2 )$ stands for the following symbol class
\begin{multline*}
S(\langle (x,\xi) \rangle^m,dx^2+d\xi^2)=\Big\{a(x,\xi;h) \in C^{\infty}(\rr_x^n \times \rr_{\xi}^n,\cc) : \\ 
\forall \alpha \in \nn^{2n}, \ \sup_{0< h \leq 1}\|\langle (x,\xi) \rangle^{-m}
\partial_{x,\xi}^{\alpha} a(x,\xi;h)\|_{L^{\infty}(\rr^{2n})} <+\infty\Big\},
\end{multline*}
with a semiclassical expansion
\begin{equation}\label{S1}\inc
P(x,\xi;h) \sim \sum_{j=0}^{+\infty}{ h^j p_j(x,\xi)}, \num
\end{equation}
where for all $j \in \nn$, $p_j$ is a symbol of the class $S(\langle (x,\xi) \rangle^m,dx^2+d\xi^2)$ independent from the semiclassical parameter $h$.

Let $z \in \cc$, we assume that there exists a function $q_0 \in C_b^{\infty}(\rr^{2n},\cc)$,
where $C_b^{\infty}(\rr^{2n},\cc)$ stands for the set of bounded complex-valued functions on $\rr^{2n}$ with all derivatives bounded, 
and a bicharacteristic curve, $t \in [a,b] \mapsto \gamma(t)$, of the real part $\textrm{Re}( q_0 (p_0-z))$ of the symbol 
$ q_0 (p_0-z)$, with $a<b$, such that  \inc
\begin{multline*}\label{S2}
\forall t \in [a,b], \ q_0\big(\gamma(t) \big) \neq 0 \textrm{ and } \\ \textrm{Im}\big[q_0(\gamma(a))\big(p_0(\gamma(a))-z\big)\big]>0>
\textrm{Im}\big[q_0(\gamma(b))\big(p_0(\gamma(b) )-z\big)\big]. \num
\end{multline*}

\bigskip

\begin{theorem}\label{psi}
Under these assumptions \emph{(\ref{S1})} and \emph{(\ref{S2})}, for all open neighbourhood~$V$ of the compact set $\gamma([a,b])$ in $\rr^{2n}$ and for all $N \in \nn$, there exist $h_0>0$ and 
$(u_h)_{0<h \leq h_0}$ a semiclassical family in $\mathcal{S}(\rr^n)$ such that  
$$
\|u_h\|_{L^2(\rr^n)}=1,\ \emph{\textrm{FS}}\big{(}(u_h)_{0<h \leq h_0}\big{)} \subset  \overline{V} \ \textsl{and } \|P(x,h\xi;h)^w u_h-z u_h\|_{L^2(\rr^n)}=O(h^N),$$
when $h \rightarrow 0^+$.
\end{theorem}

The notation $\textrm{FS}\big{(}(u_h)_{0<h \leq h_0}\big{)}$ stands for the \textit{frequency set} of the semiclassical family $(u_h)_{0<h \leq h_0}$ defined as the complement 
in $\rr^{2n}$ of the set composed by the points $(x_0,\xi_0) \in \rr^{2n}$, for which there exists a symbol $\chi_0(x,\xi;h) \in S(1,dx^2+d\xi^2)$ such that 
$$\chi_0(x_0,\xi_0;h)=1 \textrm{ and } \|\chi_0(x,h\xi;h)^w u_h\|_{L^2(\rr^n)}=O(h^{\infty}),$$
when $h \rightarrow 0^+$.

This existence result of semiclassical quasimodes is an adaptation in a semiclassical setting of the proof given by L. Hörmander in \cite{hormander} for proving 
that the condition $(\Psi)$ is a necessary condition for the solvability of a pseudodifferential 
operator (Theorem 26.4.7 in \cite{hormander}). The existence of this result has been first mentioned in \cite{dencker}. A complete proof of this adaptation in a semiclassical 
setting is given in \cite{these}. This result
shows that when the principal symbol $p_0-z$ of the symbol $P-z$ violates the condition $(\overline{\Psi})$, there exists in this point $z$ some semiclassical quasimodes 
inducing the presence of semiclassical pseudospectrum of infinite index for the semiclassical operator~$P(x,h\xi;h)^w$.

\bigskip

\noindent
\textbf{Condition \mathversion{bold}$(\overline{\Psi})$.} A complex-valued function $p \in C^{\infty}(\rr^{2n},\cc)$ fulfills the condition~$(\overline{\Psi})$ if there is no complex-valued 
function $q \in C^{\infty}(\rr^{2n},\cc)$ such that the imaginary part $\textrm{Im}(qp)$ of the function $qp$ changes sign from positive values to negative ones along 
an oriented bicharacteristic of the symbol $\textrm{Re}(qp)$ on which the function $q$ does not vanish.

\bigskip

By using the characterization given in the previous section for the interior of the numerical range $\mathring{\Sigma}(q)$ (see (\ref{48}) and (\ref{50})), we are now going to prove that the principal symbol 
$q(x,\xi)-z$
of the semiclassical operator 
$$q(x,h\xi)^w-z,$$  
violates the condition $(\overline{\Psi})$ for all $z$ in $\mathring{\Sigma}(q)$. This violation of the condition $(\overline{\Psi})$ will
induce in view of the theorem \ref{psi} that for all $z \in \mathring{\Sigma}(q)$ and $N \in \nn$, we can find a semiclassical quasimode $(u_h)_{0<h \leq h_0} \in \mathcal{S}(\rr^n)$, with $h_0>0$, 
verifying  
$$\|u_h\|_{L^2(\rr^n)}=1 \textrm{ and } \|q(x,h\xi)^w u_h -z u_h\|_{L^2(\rr^n)}=O(h^N) \textrm{ when } h \rightarrow 0^+,$$
which will end the proof of Theorem \ref{21}.

\bigskip

Let us consider $z \in \mathring{\Sigma}(q)$. We are now going to prove that there is actually a violation of the condition $(\overline{\Psi})$ for the symbol $q-z$. 
According to (\ref{48}) and (\ref{50}), there are two cases to separate.

\medskip
\noindent
\textbf{Case 1.} Let us assume that there exists $(x_0,\xi_0) \in \rr^{2n}$ such that  
\begin{equation}\label{2.3.10}\inc
z=q(x_0,\xi_0), \  \{\textrm{Re}(q-z),\textrm{Im}(q-z)\}(x_0,\xi_0)=\{\textrm{Re } q,\textrm{Im } q\}(x_0,\xi_0)<0. \num
\end{equation}
By considering the solution of the following Cauchy problem
\begin{equation}\label{2.3.12}\inc
\left\lbrace \begin{array}{l}
Y'(t)=H_{\textrm{Re } q}\big{(}Y(t)\big{)} \\
Y(0)=(x_0,\xi_0),
\end{array}\right. \num
\end{equation}
we define the following function 
\begin{equation}\label{2.3.13}\inc
f(t)=\textrm{Im } q\big{(}Y(t)\big{)}-\textrm{Im } q(x_0,\xi_0). \num
\end{equation}
As mentioned before, $(\ref{2.3.12})$ is a linear Cauchy problem. It follows that its solution $Y$ is global and that the function $f$ is well-defined on $\rr$. 
A direct computation using $(\ref{2.3.12})$ and $(\ref{2.3.13})$ gives that for all $t \in \rr$, 
\begin{equation}\label{s1}\inc
f'(t)=\{\textrm{Re } q,\textrm{Im } q\}\big{(}Y(t)\big{)}.\num
\end{equation}
Since from $(\ref{2.3.10})$, $(\ref{2.3.12})$, (\ref{2.3.13}) and (\ref{s1}),  
$$f(0)=0, \ f'(0)=\{\textrm{Re } q,\textrm{Im } q\}(x_0,\xi_0)<0$$
and $H_{\textrm{Re }q -\textrm{Re }z}=H_{\textrm{Re }q},$ we deduce in this first case that the imaginary part of the function $q-z$ changes sign, at the first order, from positive values to negative ones along 
the oriented bicharacteristic $Y$ of the symbol $\textrm{Re }q -\textrm{Re }z$. This proves that the symbol $q-z$ actually violates the condition $(\overline{\Psi})$.

\medskip
\noindent
\textbf{Case 2.} Let us now assume that there exists $(x_0,\xi_0) \in \rr^{2n}$ such that 
\begin{equation}\label{2.3.11}\inc
z=q(x_0,\xi_0), \  \{\textrm{Re}(q-z),\textrm{Im}(q-z)\}(x_0,\xi_0)=\{\textrm{Re } q,\textrm{Im } q\}(x_0,\xi_0)>0. \num
\end{equation}
We consider as in the previous case, the global solution $Y$ of the Cauchy problem (\ref{2.3.12}) and the function $f$ defined in (\ref{2.3.13}).
Since from $(\ref{2.3.12})$, (\ref{2.3.13}), $(\ref{s1})$ and (\ref{2.3.11}),  
\begin{equation}\label{s2}\inc
f(0)=0, \ f'(0)=\{\textrm{Re } q,\textrm{Im } q\}(x_0,\xi_0)>0, \num
\end{equation}
we deduce this time that the imaginary part of the function $q-z$ also changes sign, at the first order, along 
the oriented bicharacteristic $Y$ of the symbol $\textrm{Re }q -\textrm{Re }z$. Nevertheless, this change of sign is done in the \og wrong \fg \ way. It is a change of sign from negative values to positive ones,
which does not induce directly a violation of the condition~$(\overline{\Psi})$. To check that there is actually a violation of the condition~$(\overline{\Psi})$ in this second case, we need 
to study more precisely the behaviour of the function $\textrm{Im }q-\textrm{Im }z$ along this bicharacteristic $Y$.

We deduce from (\ref{s2}) that there exists $\varepsilon>0$ such that
$$\forall t \in [-\varepsilon, \varepsilon], \ f'(t)>0,$$
which induces that 
\begin{equation}\label{2.3.14}\inc
f(\varepsilon)>0 \ \textrm{and} \ f(-\varepsilon)<0, \num
\end{equation}
since from $(\ref{s2})$, $f(0)=0$.
By using the following lemma, we obtain that for all $\delta>0$, there exists a time $t_0(\delta)>\varepsilon$ such that  
\begin{equation}\label{2.3.15}\inc
|Y\big{(}t_0(\delta)\big{)}-Y(-\eps)|<\delta. \num
\end{equation}
\begin{figure}[ht]
\caption{}
\centerline{\includegraphics[scale=0.9]{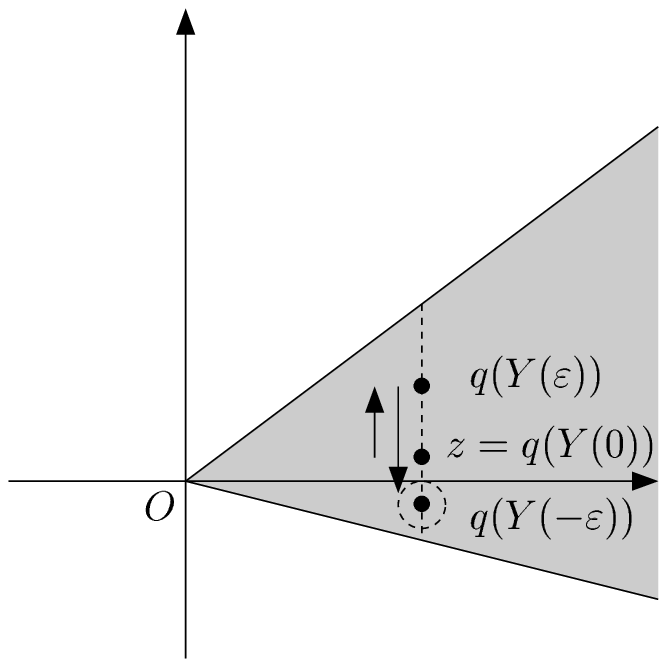}}
\end{figure}

\bigskip

\begin{lemma}\label{s3} 
If $Y(t)=(x(t),\xi(t))$ is the $C^{\infty}(\rr,\rr^{2n})$ function solving the linear system of ordinary differential equations
$$Y'(t)=H_{\emph{\textrm{Re }}q}\big{(}Y(t)\big{)},$$
where $\emph{\textrm{Re }}q$ is the symbol defined in \emph{(\ref{51})}, then we have 
\begin{multline*}
\forall t_0 \in \rr, \forall \eps>0, \forall M>0, \exists T_1>M, \exists T_2>M, \\ |Y(t_0)-Y(t_0+T_1)|<\eps \ \textrm{and} \ 
|Y(t_0)-Y(t_0-T_2)|<\eps.
\end{multline*}
\end{lemma}

\bigskip
\noindent
\textit{Proof of Lemma \ref{s3}}. If $Y(t_0)=(a_1,...,a_n,b_1,...,b_n) \in \rr^{2n}$, we deduce from (\ref{51}) that the function $Y(t)=(x(t),\xi(t))$ solves the following Cauchy problem
$$ \forall j=1,...,n, \quad
\left\lbrace \begin{array}{l}
x_j'(t)=2 \lambda_j \xi_j(t) \\
\xi_j'(t)=-2 \lambda_j x_j(t) \\
x_j(t_0)=a_j \\
\xi_j(t_0)=b_j.
\end{array}\right.$$
It follows that for all $j=1,...,n$ and $t \in \rr$, 
\begin{equation}\label{2.3.16}\inc 
\left\lbrace \begin{array}{l}
x_j(t)=b_j \sin\big{(}2(t-t_0)\lambda_j\big{)}+a_j \cos\big{(}2(t-t_0)\lambda_j\big{)} \\
\xi_j(t)=b_j \cos\big{(}2(t-t_0)\lambda_j\big{)}-a_j \sin\big{(}2(t-t_0)\lambda_j\big{)}.
\end{array}\right. \num
\end{equation}
Setting $\beta_j=\lambda_j/\pi$ for all $j=1,...,n$, we need to study two different cases.

\medskip
\noindent
\textbf{Case 1:} $\forall j \in \{1,...,n\}, \ \beta_j \in \mathbb{Q}$. In this case, the function $Y$ is periodic and the result of Lemma \ref{s3} is obvious.

\medskip
\noindent
\textbf{Case 2:} $(\beta_1,...,\beta_n) \not\in \mathbb{Q}^n$. In this second case, we use the following classical result of rational approximation: 
$\forall \eps>0$, $\forall (\theta_1,...,\theta_n) \in \rr^n \setminus \mathbb{Q}^n$, $\exists p_1,...,p_n \in \mathbb{Z}$, $\exists q \in \mathbb{N}^*$ such that 
$$0< \sup_{j=1,...,n}{\left|\theta_j-\frac{p_j}{q} \right|}<\frac{\eps}{q}.$$
If $0<\eps_1<1/2$, we can therefore find some integers $p_{1,1},...,p_{1,n} \in \mathbb{Z}$ and $q_{\eps_1} \in \mathbb{N}^*$ such that
$$0< \sup_{j=1,...,n}{\left|q_{\eps_1}\beta_j-p_{1,j} \right|}<\eps_1.$$
If  
$$\eps_2=\frac{1}{2}\sup_{j=1,...,n}{|q_{\eps_1}\beta_j-p_{1,j}|}>0,$$
using again this result of rational approximation, we can find some other integers $p_{2,1},...,p_{2,n} \in \mathbb{Z}$ and $q_{\eps_2} \in \mathbb{N}^*$ such that
$$0< \sup_{j=1,...,n}{\left|q_{\eps_2}\beta_j-p_{2,j} \right|}<\eps_2.$$
By using this process, we build some sequences $(p_{m,j})_{m \in \mathbb{N}^*}$ of $\mathbb{Z}$ for $j=1,...,n$, $(\eps_m)_{m \in
\mathbb{N}^*}$ of $\rr_+^*$ and $(q_{\eps_m})_{m \in \mathbb{N}^*}$ of $\nn^*$ such that for all $m \geq 2$,
\begin{equation}\label{2.3.17}\inc
0< \sup_{j=1,...,n}{\left|q_{\eps_m}\beta_j-p_{m,j} \right|}<\eps_m= \frac{1}{2}\sup_{j=1,...,n}{\left|q_{\eps_{m-1}}\beta_j-p_{m-1,j} 
\right|} \num
\end{equation}
and 
\begin{equation}\label{2.3.17bis}\inc
0<\eps_m < \frac{1}{2^{m-1}} \eps_1. \num
\end{equation}
The elements of the sequence $(q_{\eps_m})_{m \in \mathbb{N}^*}$ are necessary two by two different. Indeed, if $q_{\eps_k}=q_{\eps_l}$ for $k<l$, this would imply 
according to (\ref{2.3.17}) and (\ref{2.3.17bis}) that
$$\forall j=1,...,n, \ |p_{k,j}-p_{l,j}| \leq |q_{\eps_k} \beta_j- p_{k,j}|+|q_{\eps_l} \beta_j- p_{l,j}|< \eps_k+\eps_l<1,$$
because $0<\eps_1<1/2$, which would induce that $\forall j=1,...,n, \ p_{k,j}=p_{l,j}$ because $p_{k,j}$ and $p_{l,j}$ are some integers; and would contradict (\ref{2.3.17}) because 
$$0<\sup_{j=1,...,n}{\left|q_{\eps_l}\beta_j-p_{l,j} \right|}<\eps_l \leq \frac{1}{2}\sup_{j=1,...,n}{\left|q_{\eps_k}\beta_j-p_{k,j} \right|}.$$
Since the sequence $(q_{\eps_m})_{m \in \nn^*}$ is composed of integers two by two different, we can assume after a possible extraction that $q_{\eps_m} \rightarrow +\infty$ 
when $m \rightarrow +\infty$. We deduce from (\ref{2.3.16}), (\ref{2.3.17}) and (\ref{2.3.17bis}) that 
$$Y(t_0+q_{\eps_m}) \rightarrow Y(t_0) \ \textrm{when} \ m \rightarrow +\infty.$$ 
Then, considering $(\tilde{\beta_1},...,\tilde{\beta_n})=(-\beta_1,...,-\beta_n)$, we obtain by using the same method a sequence $(\tilde{q}_{\eps_m})_{m \in
\nn^*}$ of integers such that $\tilde{q}_{\eps_m} \rightarrow +\infty$ and 
$$Y(t_0-\tilde{q}_{\eps_m}) \rightarrow Y(t_0) \ \textrm{when} \ m \rightarrow
+\infty.$$ This ends the proof of Lemma \ref{s3}. $\Box$

\bigskip

Since from (\ref{2.3.14}), $f(-\eps)<0$, we deduce from (\ref{2.3.13}) and (\ref{2.3.15}) that there exists $t_0 > \eps$ such that $f(t_0)$ is arbitrarily close to $f(-\eps)$. It follows 
in particular that we can find $t_0 > \eps$ such that $f(t_0)<0$. Since from (\ref{2.3.14}), $f(\eps)>0$ and $f(t_0)<0$, we deduce from (\ref{2.3.13}) and (\ref{2.3.11}) that the function 
$$t \mapsto \textrm{Im }q\big{(}Y(t)\big{)}-\textrm{Im } z,$$ 
changes sign from positive values to negative ones on the interval $[\eps,t_0]$. This proves that the imaginary part of the function $q-z$ actually changes sign from positive values to negative 
ones along the oriented bicharacteristic $Y$ of the symbol $\textrm{Re }q -\textrm{Re }z$; and that the symbol $q-z$ also violates in this second case the condition $(\overline{\Psi})$.
This ends the proof of Theorem \ref{21}.

\subsubsection*{\ref{46.5}.c. Another proof for the existence of semiclassical quasimodes}
In the following lines, we give another proof for the existence of semiclassical quasimodes in some points of the numerical range's interior. The result proved in this section
is weaker than the one given by the theorem \ref{21}, since we prove the existence of semiclassical quasimodes in every point of the numerical range's interior without a finite number 
of particular half-lines.

Let us consider a \textit{non-normal} elliptic quadratic differential operator
\begin{equation}\label{t1}\inc
q(x,\xi)^w : B \rightarrow L^2(\rr^n), \num
\end{equation}
in dimension $n \geq 2$. We assume, as before, that (\ref{51}) is fulfilled.
Using that the quadratic form $\textrm{Re } q$ is positive definite, we can simultaneously reduce the two quadratic forms $\textrm{Re } q$ and $\textrm{Im } q$ by choosing an 
isomorphism $P$ of $\rr^{2n}$ such that in the new coordinates $y=P^{-1}(x,\xi)$, 
\begin{equation}\label{2.3.45}\inc
r_1(y)=\textrm{Re } q(Py)=\sum_{j=1}^{2n}{y_j^2}, \quad r_2(y)=\textrm{Im } q(Py)=\sum_{j=1}^{2n}{\alpha_j y_j^2}, \num
\end{equation}
with $\alpha_1 \leq ...\leq \alpha_n.$ Let us study when the differential forms $dr_1(y)$ and $dr_2(y)$ are linearly dependent on $\rr$ i.e.
when there exist $(\lambda,\mu) \in \rr^2 \setminus \{(0,0)\}$ such that
\begin{equation}\label{2.3.46}\inc
\lambda dr_1(y)+\mu d r_2(y)=0. \num
\end{equation}
It follows from (\ref{2.3.45}) and (\ref{2.3.46}) that for all $j=1,...,2n$, 
\begin{equation}\label{2.3.47}\inc
(\lambda+\mu \alpha_j)y_j=0. \num
\end{equation}
If $y \neq 0$, then there exists $j_0 \in \{1,...,2n\}$ such that $y_{j_0} \neq 0$. This implies that  
\begin{equation}\label{2.3.48}\inc
\lambda+\mu \alpha_{j_0}=0. \num
\end{equation}
We deduce from (\ref{2.3.47}) and (\ref{2.3.48}) that $y_j=0$ if $\alpha_j \neq \alpha_{j_0}$. Thus, we obtain that if 
$$z \in \mathring{\Sigma}(q) \setminus \big{(}(1+i \alpha_1)\rr_+^* \cup ... \cup (1+i \alpha_n)\rr_+^* \big{)},$$
then the differential forms $d\textrm{Re }q$ and $d\textrm{Im}q$ are linearly independent on $\rr$ in every point of the set $q^{-1}(z)$. 
\begin{figure}[ht]
\caption{}
\centerline{\includegraphics[scale=0.9]{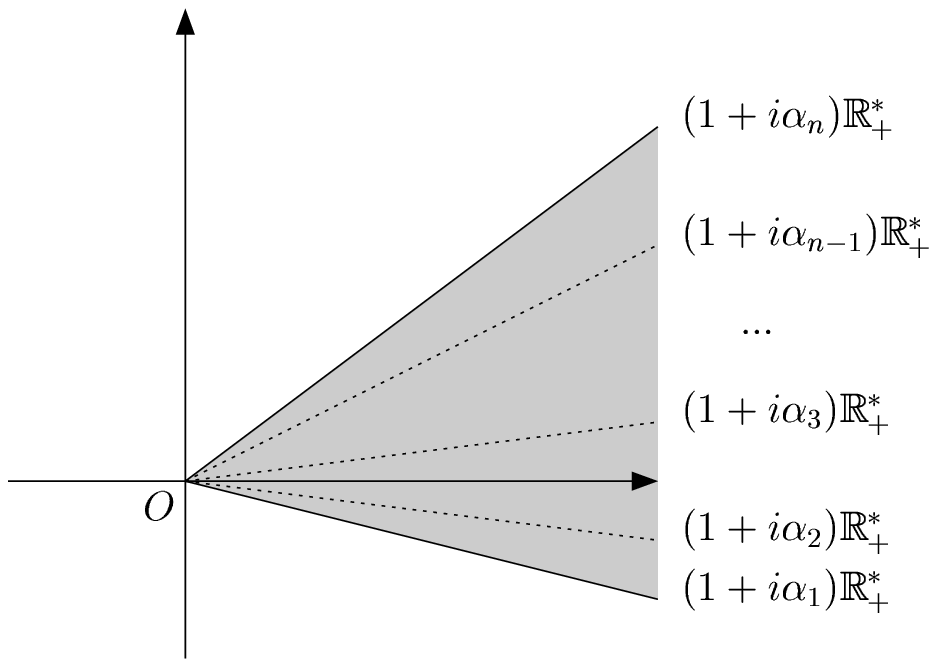}}
\end{figure}

Let us consider such a point
$$z \in \mathring{\Sigma}(q) \setminus \big{(}(1+i \alpha_1)\rr_+^* \cup ... \cup (1+i \alpha_n)\rr_+^* \big{)}.$$
Since the dimension $n \geq 2$, we can apply the lemma $3.1$ in \cite{dencker} (see also the lemma 8.1 in \cite{melin}). It follows that for any compact, 
connected component $\Gamma$ of $q^{-1}(z)$, we have
\begin{equation}\label{2.3.49}\inc
\int_{\Gamma}{\{\textrm{Re } q, \textrm{Im } q\}(\rho) \lambda_{q,z}(d\rho)}=0, \num
\end{equation}
where $\lambda_{q,z}$ stands for the Liouville measure on $q^{-1}(z)$, 
$$\lambda_{q,z} \wedge d\textrm{Re }q \wedge d\textrm{Im }q=\frac{\sigma^n}{n!}.$$
The set $q^{-1}(z)$ is a non-empty submanifold of codimension $2$ in $\rr^{2n}$. We deduce from (\ref{48}) and (\ref{50}) that there exist $(x_0,\xi_0) \in q^{-1}(z)$ such that 
\begin{equation}\label{2.3.50}\inc
\{\textrm{Re } q, \textrm{Im } q\}(x_0,\xi_0) \neq 0. \num
\end{equation}
Then, it follows from (\ref{2.3.49}) and (\ref{2.3.50}) that there necessary exists $(\tilde{x}_0,\tilde{\xi}_0) \in q^{-1}(z)$ such that 
\begin{equation}\label{t10}\inc
\{\textrm{Re } q, \textrm{Im } q\}(\tilde{x}_0,\tilde{\xi}_0) < 0.\num
\end{equation} 
Under this condition (\ref{t10}), we can use the reasoning given in the first studied case (see (\ref{2.3.10})) to prove that the imaginary part of the function $q-z$ changes sign, at the first order, 
from positive values to negative ones along an oriented bicharacteristic of the symbol $\textrm{Re }q -\textrm{Re }z$. This induces that the symbol $q-z$ violates 
the condition~$(\overline{\Psi})$; and we can conclude by using the theorem \ref{psi}. Let us mention that we can also directly use the existence result of semiclassical 
quasimodes given by M. Zworski in \cite{zworski1} and \cite{zworski2}. This second proof gives the existence of semiclassical quasimodes in every point belonging to the set  
$$\mathring{\Sigma}(q) \setminus \big{(}(1+i \alpha_1)\rr_+^* \cup ... \cup (1+i \alpha_n)\rr_+^* \big{)}.$$

\subsubsection{On the pseudospectrum at the boundary of the numerical range}

In this section, we give a proof  of the theorem \ref{27}. Let us consider a \textit{non-normal} elliptic quadratic differential operator
$$q(x,\xi)^w : B \rightarrow L^2(\rr^n),$$
in dimension $n \geq 1$. We assume that $\Sigma(q) \neq \cc$, and that its Weyl symbol $q(x,\xi)$ is of \textit{finite} order $k_j$ on a half-line $\Delta_j$, $j \in \{1,2\}$
(See the definition given in (\ref{t10.5})), which composes 
the boundary of its numerical range
\begin{equation}\label{t11}\inc
\partial \Sigma(q)=\{0\}  \sqcup \Delta_1 \sqcup \Delta_2. \num
\end{equation} 
As we have already done several times, we can reduce our study to case where (\ref{51}) is fulfilled.

\bigskip

\noindent
\textit{Proof of Theorem \ref{27}.}
Let us consider the following symbol belonging to the $C_b^{\infty}(\rr^{2n},\cc)$ space, composed of bounded complex-valued functions on $\rr^{2n}$ with all derivatives bounded
\begin{equation}\label{2.3.62}\inc 
r(x,\xi)=\frac{q(x,\xi)-z}{1+x^2+\xi^2}, \num
\end{equation}
with $z \in \Delta_j$. Setting 
$\tilde{\Sigma}(r)=\overline{r(\rr^{2n})},$
we can first notice that 
$$z \in \partial \Sigma(q) \setminus \{0\} \Rightarrow 0 \in \partial \tilde{\Sigma}(r).$$ 
Let us also notice that the symbol $r$ fulfills the principal-type condition in $0$. Indeed, if $(x_0,\xi_0) \in \rr^{2n}$ was such that $r(x_0,\xi_0)=0$ and $dr(x_0,\xi_0)=0$, 
we would get from (\ref{2.3.62}) that  
\begin{equation}\label{2.3.63}\inc
dq(x_0,\xi_0)=0. \num
\end{equation}
Since from (\ref{51}) and (\ref{2.3.63}), we have
$$d\textrm{Re } q(x_0,\xi_0)=2\sum_{j=1}^{n}{\lambda_j \big{(}(x_0)_j dx_j+ (\xi_0)_j d\xi_j \big{)}}=0,$$
this would imply that 
$$(x_0,\xi_0)=(0,0), \ q(x_0,\xi_0)=0,$$
because $q$ is a quadratic form and that $\lambda_j>0$ for all $j=1,...,n$. On the other hand, since $r(x_0,\xi_0)=0$, we get from (\ref{2.3.62}) that $q(x_0,\xi_0)=z \neq 0$ because  
$$z \in \Delta_j \subset \partial\Sigma(q) \setminus \{0\},$$ 
which induces a contradiction. It follows that the symbol $r$ actually fulfills the principal-type condition in $0$.
Let us notice that, since symbol $q$ is of finite order~$k_j$ in $z$, this induces in view of (\ref{2.3.62}) that the symbol $r$ is also of finite order $k_j$ in $0$.
On the other hand, we deduce from (\ref{51}) and (\ref{2.3.62}) that the set
$$\{(x,\xi) \in \rr^{2n} : r(x,\xi)=0\}=\{(x,\xi) \in \rr^{2n} : q(x,\xi)=z\},$$
is compact. Under these conditions, we can apply the theorem 1.4 in \cite{dencker}, which proves that the integer~$k_j$ is \textit{even} and gives the existence of positive constants $h_0$ and $C_1$ 
such that 
\begin{equation}\label{2.3.64}\inc
\forall \ 0<h<h_0, \forall u \in \mathcal{S}(\rr^n), \ \|r(x,h\xi)^w u\|_{\lde} \geq C_1 h^{\frac{k_j}{k_j+1}}\|u\|_{\lde}. \num
\end{equation}

\bigskip
\noindent
\textit{Remark.} We did not check the dynamical condition $(1.7)$ in~\cite{dencker}, because this assumption is not necessary
for the proof of Theorem~1.4. Indeed, this proof only use a part of the proof of lemma 4.1 in \cite{dencker} (a part of the second paragraph), where this condition (1.7) is not needed.

\bigskip
\noindent
By using some results of symbolic calculus given by Theorem 18.5.4 in \cite{hormander} and (\ref{2.3.62}), we can write
\begin{equation}\label{2.3.65}\inc
r(x,h\xi)^w(1+x^2+h^2\xi^2)^w=q(x,h\xi)^w-z + h r_1(x,h\xi)^w +h^2 r_2(x,h\xi)^w, \num
\end{equation}
with 
\begin{equation}\label{2.3.66}\inc
r_1(x,\xi)=-i x \frac{\partial r}{\partial \xi}(x,\xi)+i \xi \frac{\partial r}{\partial x}(x,\xi) \num
\end{equation}
and
\begin{equation}\label{2.3.67}\inc
r_2(x,\xi)=-\frac{1}{2}\frac{\partial^2 r}{\partial x^2}(x,\xi)-\frac{1}{2}\frac{\partial^2 r}{\partial \xi^2}(x,\xi). \num
\end{equation}
We can easily check from (\ref{2.3.62}) that these functions $r_1$ and $r_2$ belong to the space $C_b^{\infty}(\rr^{2n},\cc)$, and we deduce from the Calder\'on-Vaillancourt theorem 
that there exists a positive constant $C_2$ such that for all $u \in \mathcal{S}(\rr^n)$ and $0<h \leq 1$,
\begin{equation}\label{2.3.68}\inc
\|r_1(x,h\xi)^w u\|_{L^2} \leq C_2 \|u\|_{L^2} \textrm{ and }  \|r_2(x,h\xi)^w u\|_{L^2} \leq C_2 \|u\|_{L^2}. \num
\end{equation}
It follows from (\ref{2.3.64}), (\ref{2.3.65}), (\ref{2.3.68}) and the triangular inequality that for all $u \in \mathcal{S}(\rr^n)$ and $0<h<h_0$,
\begin{align*}
& \ C_1 h^{\frac{k_j}{k_j+1}} \|(1+x^2+h^2 \xi^2)^w u\|_{L^2(\rr^n)} \\
\leq & \ \|r(x,h\xi)^w (1+x^2+h^2 \xi^2)^w u \|_{L^2(\rr^n)} \\
\leq & \ \|q(x,h\xi)^w u -z u \|_{L^2(\rr^n)} + C_2 h(1+h)\|u \|_{L^2(\rr^n)}. 
\end{align*}
Since from the Cauchy-Schwarz inequality, we have for all $u \in \mathcal{S}(\rr^n)$ and $0<h \leq 1$,
\begin{align*}
\|u \|_{L^2(\rr^n)}^2 \leq & \ \|u \|_{L^2(\rr^n)}^2+\|xu \|_{L^2(\rr^n)}^2+\|hD_xu \|_{L^2(\rr^n)}^2\\
= & \ \big((1+x^2+h^2 \xi^2)^w u,u\big)_{L^2(\rr^n)} \\ 
\leq & \ \|(1+x^2+h^2 \xi^2)^w u \|_{L^2(\rr^n)}\|u \|_{L^2(\rr^n)}, 
\end{align*}
we obtain that for all $u \in \mathcal{S}(\rr^n)$ and $0<h <h_0$,
\begin{equation}\label{2.3.72}\inc
C_1 h^{\frac{k_j}{k_j+1}} \|u\|_{L^2(\rr^n)} \leq \|q(x,h\xi)^w u -z u \|_{L^2(\rr^n)} + C_2 h(1+h)\|u \|_{L^2(\rr^n)}. \num
\end{equation}
Since $k_j \geq 1$, we deduce from (\ref{2.3.72}) that there exist some positive constants $h'_0$ and $C_{3}$ such that 
for all $0<h<h'_0$ and $u \in \mathcal{S}(\rr^n)$,  
$$\|q(x,h\xi)^w u - z u\|_{\lde} \geq C_{3} h^{\frac{k_j}{k_j+1}}\|u\|_{\lde}.$$
Using that the Schwartz space $\mathcal{S}(\rr^n)$ is dense in $B$ and that the operator 
$$q(x,h\xi)^w+z,$$
is a Fredholm operator of index $0$, we obtain that for all $0<h<h_0'$,
$$\big\|\big(q(x,h\xi)^w-z\big)^{-1}\big\| \leq C_3^{-1}h^{-\frac{k_j}{k_j+1}},$$
which ends the proof of Theorem \ref{27}. $\Box$

\bigskip

About the case of \textit{infinite} order, the situation is much more complicated. As mentioned before, we cannot expect to prove a stronger result than an 
absence of semiclassical pseudospectrum of index 1, but we can actually prove that there is never some semiclassical pseudospectrum of index 1 on 
every half-line of infinite order, by using a result of exponential decay in time for the norm of contraction semigroups generated 
by elliptic quadratic differential operators proved in \cite{karel6}.

The result proved in \cite{karel6} shows that the norm of a contraction semigroup 
$$\|e^{tq(x,\xi)^w}\|_{\mathcal{L}(L^2)}, \ t \geq 0,$$
generated by an elliptic quadratic differential operator $q(x,\xi)^w$ with a Weyl symbol verifying
$$\textrm{Re }q \leq 0, \ \exists (x_0,\xi_0) \in \rr^{2n}, \ \textrm{Re }q(x_0,\xi_0) \neq 0,$$
decreases exponentially in time
\begin{equation}\label{stefy1}\inc
\exists M,a>0, \forall t \geq 0, \ \|e^{tq(x,\xi)^w}\|_{\mathcal{L}(L^2)} \leq M e^{-a t}.\num
\end{equation}

Let us consider a \textit{non-normal} elliptic quadratic differential operator
$$q(x,\xi)^w : B \rightarrow L^2(\rr^n),$$
in dimension $n \geq 1$ such that $\Sigma(q) \neq \cc$. 
We explain in the following lines how (\ref{stefy1}) allows to prove that there is never some semiclassical pseudospectrum of index 1 on any open half-lines composing 
the boundary of the numerical range $\partial \Sigma(q) \setminus \{0\}$.

Let $z \in \partial \Sigma(q) \setminus \{0\}$. Since the numerical range $\Sigma(q)$ is a closed angular sector with a top in $0$ and a positive opening strictly lower than $\pi$, we can find
$\eps \in \{\pm 1\}$ such that 
\begin{equation}\label{stefy2}\inc
\textrm{Re}(\eps i z^{-1}q) \leq 0, \ \exists (x_0,\xi_0) \in \rr^{2n}, \ \textrm{Re}(\eps i z^{-1}q)(x_0,\xi_0) \neq 0. \num
\end{equation} 
Using the theorem 2.8 in \cite{davies}, we obtain that for all $\eta \in \rr$,
\begin{align*}\label{2.3.85}\inc
\big{(}q(x,\xi)^w-\eta z \big{)}^{-1}= & \ -i z^{-1} \eps \big{(}\eps i \eta-\eps i z^{-1} q(x,\xi)^w \big{)}^{-1}\\
= & \ -i z^{-1} \eps \int_{0}^{+\infty}{e^{-i \eps \eta s }e^{s \eps i z^{-1} q(x,\xi)^w}}ds. \num
\end{align*}
It follows from (\ref{stefy1}) and (\ref{stefy2}) that for all $\eta \in \rr$,
\begin{align*}
\big\|\big{(}q(x,\xi)^w-\eta z\big{)}^{-1}\big\| \leq & \ |z|^{-1} \int_{0}^{+\infty}{\|e^{s \eps i z^{-1} q(x,\xi)^w}\|_{\mathcal{L}(L^2)} ds} \\
\leq & \  |z|^{-1} \int_{0}^{+\infty}{M e^{-a s}ds}= |z|^{-1} \frac{M}{a} < +\infty, 
\end{align*}
which proves the absence of semiclassical pseudospectrum of index 1 on the half-line $z \rr_+^*$. 
We can actually use the theorem 2.8 in \cite{davies} because 
$$i\rr \subset \cc \setminus \sigma\big( \eps i z^{-1} q(x,\xi)^w\big).$$
Indeed, if it was not the case, we would deduce from (\ref{15}) that there exists $u_0 \in B \setminus \{0\}$ and $\lambda_0 \in \rr$ such that 
$$\eps i z^{-1} q(x,\xi)^w u_0=i \lambda_0 u_0.$$ 
Since from (\ref{stefy2}), the quadratic form $-\textrm{Re}(\eps i z^{-1}q)$ is non-negative, we deduce from the symplectic invariance of the Weyl quantization and the theorem 21.5.3 in \cite{hormander} 
that there exists a metaplectic operator $U$ such that  
\begin{equation}\label{stefy3}\inc
-\big[\textrm{Re}\big(\eps i z^{-1}q(x,\xi)\big)\big]^w=U^{-1} \Big(\sum_{j=1}^{k}{\lambda_j(D_{x_j}^2+x_j^2)}
+\sum_{j=k+1}^{k+l}{x_j^2} \Big) U, \num
\end{equation}
with $k,l \in \nn$ and $\lambda_j>0$ for all $j=1,...,k$. By using that $U$ is a unitary operator on $L^2(\rr^{n})$, we obtain that 
\begin{align*}
0= & \ -\textrm{Re}(i \lambda_0 u_0,u_0)_{L^2}\\
= & \ -\textrm{Re}\big(\eps i z^{-1}q(x,\xi)^w u_0,u_0\big)_{L^2}\\
= & \ -\big{(} \big[\textrm{Re}\big{(}\eps i z^{-1}q(x,\xi)\big{)}\big]^w u_0,u_0\big{)}_{L^2} \\
= & \ \sum_{j=1}^{k}{\lambda_j \big{(}\|D_{x_j} U u_0\|_{L^2}^2+\|x_j U u_0\|_{L^2}^2 \big{)}}+\sum_{j=k+1}^{k+l}{\|x_j U u_0\|_{L^2}^2},
\end{align*}
which induces that $u_0=0$, because from (\ref{stefy2}) and (\ref{stefy3}), $k+l \geq 1$. It follows from (\ref{15}) that there exists $\eps_0>0$ such that
$$\sigma\big(\eps i z^{-1} q(x,\xi)^w\big) \subset \{z \in \cc : \textrm{Re } z \leq -\eps_0\}.$$

\bigskip
\bigskip

\noindent
\textsc{Department of Mathematics, University of California, Evans Hall,
Berkeley, CA 94720, USA}\\
\textit{E-mail address:} \textbf{karel@math.berkeley.edu}

\end{document}